%% file: main.tex
\documentclass{amsart}
\input{Macro.tex}

\usepackage{natbib}
\usepackage{graphicx}
\usepackage{xcolor}
\usepackage{amsmath, amsfonts, amssymb, amsthm, color, dsfont}
\usepackage{enumitem}
\usepackage[utf8]{inputenc}
\usepackage[colorlinks,linkcolor=blue,citecolor=blue]{hyperref}
\usepackage[all, knot, color]{xy}
\xyoption{arc}
\xyoption{line}
\usepackage{booktabs}
\UseCrayolaColors
\usepackage{tikz}
\usetikzlibrary{intersections}
\usetikzlibrary{cd}
\tikzset{
  dotted stroke/.style={thick, densely dotted}
}

\newtheorem{propdef}[thm]{Proposition--Definition}

\title{Categorical orbifold theory}

\author{Chongying Dong}
\address{Department of Mathematics, University of California, Santa Cruz, CA 95064}
\email{dong@ucsc.edu}

\author{Siu-Hung Ng}
\address{Department of Mathematics, Louisiana State University, Baton Rouge, LA 70817}
\email{rng@math.lsu.edu}

\author{Li Ren}
\address{Department of Mathematics, Sichuan University, Chengdu, 610064, China}
\email{RenL@scu.edu.cn}

\author{Feng Xu}
\address{Department of Mathematics, University of California, Riverside, CA, 92521}
\email{xufeng@math.ucr.edu}

\begin{document}

\begin{abstract}
We investigate the orbifold theory for condensable algebras in fusion
categories $\EuScript{C}$. We prove that the algebra automorphism group of a condensable
algebra $A$ is finite and explicitly characterize its elements in terms
of the fusion action and indecomposable
$\EuScript{C}_{\!A}^0$-bimodule decomposition of $\EuScript{C}_{\!A}$.
We further show that a finite group action on $A$ gives rise to a group
extension of the category $\EuScript{C}_{\!A}^0$ in $\EuScript{C}_{\!A}$.
Conversely, every such group extension arises from a finite group of
algebra automorphisms of $A$. Finally, we also construct a canonical
maximal group extension of $\EuScript{C}_{\!A}^0$, thereby providing a
categorical framework for the classical correspondence among finite
group actions, fixed-point algebras, and twisted modules.
\end{abstract}
\maketitle

\section{Introduction}

Orbifold theory provides a fundamental connection between finite group actions, fixed-point algebras, and twisted representations. In the theory of vertex operator algebras, when a finite group $G$ acts on a vertex operator algebra $V$, the fixed-point subalgebra $V^G$ is a central object of study, while $g$-twisted $V$-modules, for $g\in G$, play an essential role in understanding the representation theory of $V^G$. Under suitable finiteness and rationality assumptions, the categories of twisted $V$-modules assemble into a $G$-graded extension of the category of ordinary $V$-modules. This structure naturally suggests that the fundamental ingredients of orbifold theory should admit an intrinsic formulation in the language of tensor categories, independent of their realization through vertex operator algebras \cite{Ki2002, Ki3, McRae2021}.

The purpose of this paper is to further develop a categorical orbifold theory for condensable algebras in modular tensor categories. Let $\CC$ be a modular tensor category, and let $A$ be a condensable algebra in $\CC$. Denote by $\CA$ the fusion category of right $A$-modules and by $\CAo$ the modular tensor category of local $A$-modules. The category $\CA$ acts on the algebra object $A$ through the fusion action introduced in our previous work \cite{DNR2025}. The associated categorical Schur--Weyl duality identifies the effective fusion algebra action with $\operatorname{End}_{\CC}(A)$, while the corresponding Galois correspondence relates fusion subcategories of $\mathcal{C}_A$ containing $\mathcal{C}_A^0$ to condensable subalgebras of $A$.  Moreover, we strengthen the Galois correspondence established in \cite{DNR2025}
by removing the pseudounitarity condition. This improved result is proved and stated as Theorem  \ref{t:Galois-fusion}. These results provide the starting point for the present work. 

We first study algebra automorphisms of condensable algebras. A basic result of this paper is that the algebra automorphism group
$
\operatorname{Aut}_{\CC}^{\mathrm{alg}}(A)$
is always finite. More precisely, distinct algebra automorphisms of $A$ are linearly independent in $\operatorname{End}_{\CC}(A)$, and consequently
$$
\left|\operatorname{Aut}_{\CC}^{\mathrm{alg}}(A)\right|
\leq \dim \operatorname{End}_{\CC}(A).
$$
This strengthens the result of \cite{Ki2002} that every  compact subgroup of $\operatorname{Aut}_{\CC}^{\mathrm{alg}}(A)$ is a finite group. Our finiteness theorem shows that finite group actions on condensable algebras arise naturally from the underlying categorical structure. In fact, $\operatorname{Aut}_{\CC}^{\mathrm{alg}}(A)$ is determined completely by the $\CAo$-module structure of $\CA.$

For an algebra automorphism $g$ of $A$, there is the notion of a $g$-twisted $A$-module \cite{Ki2002}. The defining condition is a categorical analogue of the monodromy relation satisfied by twisted modules in vertex operator algebra theory. We prove that the category $\CA(g)$ of $g$-twisted $A$-modules is nontrivial for every algebra automorphism $g$. Moreover, the fusion action of every object in $\CA(g)$ determines the automorphism $g$.

Let $G$ be a subgroup of $\operatorname{Aut}_{\CC}^{\mathrm{alg}}(A)$, and let $B = A^G$ be the corresponding fixed-point subalgebra. We prove that the twisted module categories assemble into a $G$-graded fusion category
$$
\bigoplus_{g \in G} \CA(g),
$$
whose trivial component is $\CAo$. More precisely, this category coincides with the category of $A$-modules internal to the category of local $B$-modules and is therefore a $G$-extension of $\CAo$. We also recover $G$ as the full algebra automorphism group of $A$ inside the local $B$-module category.

Conversely, suppose that a fusion subcategory of $\CA$ is a $G$-extension of $\CAo$. Using the fusion action, we construct from each homogeneous component of the extension an algebra automorphism of $A$ and prove that its objects are precisely the corresponding twisted $A$-modules. Consequently, group extensions of $\CAo$ inside $\CA$ are exactly the fusion subcategories consisting of twisted $A$-modules. This establishes a categorical correspondence between finite group actions on condensable algebras and group extensions of their local module categories. 

Finally, we study group extensions of fusion subcategories from an intrinsic fusion-categorical viewpoint. Given a fusion category $\DD$ and a fusion subcategory $\DD_1$ admitting an appropriate central structure, we consider the decomposition of $\DD$ into indecomposable $\DD_1$-bimodule categories. We associate canonical elements of the Grothendieck algebra to these components and characterize the invertible ones. These elements form a finite group, and the corresponding components generate a canonical maximal group extension of $\DD_1$ inside $\DD$.

The general framework for group extensions of fusion categories was developed by Etingof, Nikshych, and Ostrik in their foundational work \cite{ENO2010}.
They showed that, for a $G$-extension $\DD=\bigoplus_{g\in G}\DD_g$ of a fusion category $\DD_1$, every homogeneous component $\DD_g$ is an invertible $\DD_1$-bimodule category. More generally, such extensions are classified  by morphisms of categorical $2$-groups $G\to\underline{\underline{\operatorname{BrPic}}}(\mathcal D_1)$, together with the corresponding coherence and obstruction data. Our categorical orbifold construction establishes a direct relationship
between the automorphism group $\Auta(A)$ of a condensable algebra $A$ 
and the group extensions of $\CAo.$ Namely,  any subgroup $G\leq\Auta(A)$ gives rise to a $G$-graded extension of $\CAo$ inside $\CA,$
whose $g$-component consists of $g$-twisted $A$-modules. Conversely, for 
every such  group extension inside $\CA$, the grading group is realized as a subgroup of $\Auta(A)$. The canonical maximal extension constructed below therefore assembles all twisted sectors arising from algebra automorphisms of $A$.

Taken together, these results provide a categorical framework for the basic structures of orbifold theory. Finite group actions, fixed-point algebras, twisted modules, and group extensions of tensor categories arise as different manifestations of the same categorical structure. The resulting theory extends classical orbifold constructions and provides a foundation for further applications to vertex operator algebras, conformal field theory, and related areas of tensor category theory.\\

{\bf Acknowledgement:}  S.-H. Ng was partially supported by the Simons Foundation MPS-TSM-00008039. L. Ren was supported by NSFC Grant 12571029.

\section{Basic setup and fusion actions}

In this section, we discuss condensable algebras and their
module categories, and recall the fusion action associated with a
condensable algebra. We refer the reader to \cite{EGNO} for general
background on fusion categories and to \cite{DNRX2025} for the construction
of fusion actions.

Let $\eM$ be a semisimple $\BC$-linear category. We denote by $\irr(\eM)$ a
complete set of nonisomorphic simple objects of $\eM$. Following the
convention in \cite{MacLane}, the morphism space of $\eM$ from $x$ to $y$
will be denoted by $\eM(x,y)$, and $[x,y]_\eM:=\dim_\BC \eM(x,y)$.

A fusion category in this paper means a semisimple $\BC$-linear rigid finite
tensor category with finite-dimensional hom-spaces. For a fusion category
$\eM$, we denote by $K_0(\eM)$ its Grothendieck ring and set
$K(\eM):=K_0(\eM)\o_\BZ \BC$. We identify the canonical basis of
$K_0(\eM)$ with $\irr(\eM)$. It is well known that $K(\eM)$ is a
semisimple associative $\BC$-algebra. If $\eM$ is pivotal,  we may assume $\eM$ is strict pivotal (cf. \cite[Thm. 2.2]{NSind}). In this paper, we only consider spherical pivotal fusion categories $\eM$, and the corresponding categorical (pivotal) dimension of any object $x$ in $\eM$ will be denoted $d(x)$. The Frobenius-Perron dimension of $x$, denoted by $d_+(x)$, is independent of the pivotal structure of $\eM$.

We will use the same graphical conventions detailed in \cite{DNRX2025}. In particular,  the evaluation $\ev_y$ and coevaluation $\coev_y$ morphisms for the left dual of $y \in \eM$ are depicted by the diagrams:
\[
\coev_y := \vcenter{\hbox{%
\begin{tikzpicture}[x=1pc, y=1pc, thick]
  \draw (0,2.1) arc[start angle=180, end angle=0, radius=1pc];
  \node[inner sep=0pt, outer sep=0pt] at (0,1.8) {$\scriptstyle y$};
  \node[inner sep=0pt, outer sep=0pt] at (2.2,1.8) {$\scriptstyle y^*$};
\end{tikzpicture}%
}}\quad, \quad
\ev_y:=\vcenter{\hbox{%
\begin{tikzpicture}[x=1pc, y=1pc, thick]
  \draw (8,2.7) arc[start angle=180, end angle=360, radius=1pc];
  \node[inner sep=0pt, outer sep=0pt] at (8,3.3) {$\scriptstyle y^*$};
  \node[inner sep=0pt, outer sep=0pt] at (10,3.3) {$\scriptstyle y$};
\end{tikzpicture}%
}}\,.%
\]
The right duality morphisms $\tilde \ev_y$ and $\tilde \coev_y$ use the same diagrams except the labeling locations of $y$ and $y^*$ are interchanged.

Let $\CC$ be a modular tensor category. We write $R$ for its braiding,
$\theta$ for its twist, and $d$ for its categorical dimension function.
For $x,y\in \irr(\CC)$, the $S$-matrix is defined by
$S^\CC_{x,y}:=\tr(R_{y,x^*}R_{x^*,y})$, where $\tr_\CC$ is the categorical (pivotal) trace of $\CC$. If $\CC$ is pseudounitary, we
always use its canonical spherical pivotal structure. Throughout this paper, we may assume $\CC$ is strict pivotal. The braidings of $\CC$ are depicted by
\[
R_{x,y}=\vcenter{\hbox{%
\begin{tikzpicture}[x=1pc, y=1pc, thick]
  \begin{scope}[even odd rule]
    \clip[overlay] (-2,-2) rectangle (2,2) (0,0) circle (2.5pt);
    \draw (0.5,1)
      .. controls (0.5,0.2) and (-0.5,-0.2) .. (-0.5,-1);
  \end{scope}
  \draw (-0.5,1)
    .. controls (-0.5,0.2) and (0.5,-0.2) .. (0.5,-1);
  \node[anchor=south, inner sep=2pt] at (-0.5,1) {$\scriptstyle x$};
  \node[anchor=south, inner sep=2pt] at (0.5,1) {$\scriptstyle y$};
\end{tikzpicture}%
}}, \qquad
R_{x,y}^{-1}=\vcenter{\hbox{%
\begin{tikzpicture}[x=1pc, y=1pc, thick]
  \begin{scope}[even odd rule]
    \clip[overlay] (-2,-2) rectangle (2,2) (0,0) circle (2.5pt);
    \draw (-0.5,1)
      .. controls (-0.5,0.2) and (0.5,-0.2) .. (0.5,-1);
  \end{scope}
  \draw (0.5,1)
    .. controls (0.5,0.2) and (-0.5,-0.2) .. (-0.5,-1);
  \node[anchor=south, inner sep=2pt] at (-0.5,1) {$\scriptstyle y$};
  \node[anchor=south, inner sep=2pt] at (0.5,1) {$\scriptstyle x$};
\end{tikzpicture}%
}}\,.
\]

We recall the type of algebra object that will be used throughout the
paper. An object $A\in \CC$ is an algebra if it is equipped with morphisms
$m_A:A\o A\to A$ and $u_A:\1\to A$ satisfying the usual associativity and
unit axioms. The algebra $A$ is called connected if $[\1,A]_\CC=1$, and
commutative if $m_AR_{A,A}=m_A$. We say that $A$ is
\emph{condensable} if $A$ is commutative, connected, $d(A)\ne 0$,
$\theta_A=\id_A$, and the pairing $\e_A m_A:A\o A\to \1$ is nondegenerate,
where $\e_A\in \CC(A,\1)$  is the retraction of $u_A$, i.e. $\e_A u_A = \id_\1$ (cf. \cite[p. 193]{KO02}). A condensable algebra in $\CC$ is called \emph{rigid} in \cite{KO02}.

 The simplicity of $A$ implies the pairing $\e_A m_A: A \o A \to \1$ is nondegenerate \cite[Lem. 1.20]{KO02}, and so one can choose $\ev_A = \e_A m_A$ (cf. \cite[p.195]{KO02}). The associativity of $m_A$ implies that $A$ is a Frobenius algebra. We will use the same graphical conventions detailed in \cite{DNRX2025}. In particular, the diagrams for the unit map $u_A$, its retraction $\varepsilon_A$, and multiplication $m_A$ of $A$ are depicted by
 \[\displaystyle
u_A:= \,\vcenter{\hbox{%
\begin{tikzpicture}[yscale=.7, x=1pc, y=1pc]
  \node[inner sep=0pt, outer sep=0pt] (unit) at (0,1.1) {$\circ$};
  \draw[dotted stroke] (unit.south) -- (0,-1.3);
\end{tikzpicture}%
}} \qquad,\qquad
\e_A:= \,\vcenter{\hbox{%
\begin{tikzpicture}[yscale=.7, x=1pc, y=1pc]
  \node[inner sep=0pt, outer sep=0pt] (counit) at (0,-1.3) {$\circ$};
  \draw[dotted stroke] (counit.north) -- (0,1.1);
\end{tikzpicture}%
}} \qquad\text{and}\qquad
m_A:= \vcenter{\hbox{%
\begin{tikzpicture}[yscale=.7, x=1pc, y=1pc]
  \draw[dotted stroke] (0,0) -- (-1,1.5);
  \draw[dotted stroke] (0,0) -- (1,1.5);
  \draw[dotted stroke] (0,0) -- (0,-1.5);
\end{tikzpicture}%
}} 
 \quad,
\]
where the dotted line will always denote the object $A$ in $\CC$. The evaluation and coevaluation of $A$ are given by
\[
  \ev_A := \vcenter{\hbox{%
\begin{tikzpicture}[yscale=.7, x=1pc, y=1pc]
  \draw[dotted stroke] (0,0) -- (-1,1.5);
  \draw[dotted stroke] (0,0) -- (1,1.5);
  \draw[dotted stroke] (0,0) -- (0,-1.5);
  \node[inner sep=0pt, outer sep=0pt] (unit) at (0,-1.6) {$\circ$};
\end{tikzpicture}%
}}   =
\vcenter{\hbox{\begin{tikzpicture}[x=1pc, y=1pc, thick]
  \draw[dotted stroke]  (8,2.7) arc[start angle=180, end angle=360, radius=1pc];
\end{tikzpicture}%
}}\quad \text{and}\quad
\coev_A = \vcenter{\hbox{\begin{tikzpicture}[x=1pc, y=1pc, thick]
  \draw[dotted stroke] (0,2.1) arc[start angle=180, end angle=0, radius=1pc];
\end{tikzpicture}%
}}\quad.
\]

Assume from now on that $A$ is a condensable algebra in $\CC$. Let $\CA$
denote the category of right $A$-modules in $\CC$. By \cite{KO02}, $\CA$
is a spherical fusion category with unit object $A$ and tensor product
$\oa$ over $A$. For $X\in \CA$, its dimension in $\CA$ is denoted by
$d_A(X)$, and $d_A(X)=d(X)/d(A)$. When the context is clear, we denote $A$
as the unit object $\1$ of $\CA$. The right $A$-module action $\mu_X$ and its chosen left $A$-action $\mu_X^+$ are depicted by the diagrams: 
\[
\mu_X=\vcenter{\hbox{%
\begin{tikzpicture}[x=1pc, y=1pc, thick]
  \draw (0,1.5) -- (0,-1.5);
  \draw[dotted stroke] (0,0) -- (1.5,1.5);
  \node[anchor=south, inner sep=2pt] at (0,1.5) {$\scriptstyle X$};
  \node[anchor=north, inner sep=2pt] at (0,-1.5) {$\scriptstyle X$};
\end{tikzpicture}%
}} \qquad\text{and}\qquad
\mu_X^+=\vcenter{\hbox{%
\begin{tikzpicture}[x=1pc, y=1pc, thick]
  \path[name path=action vertical] (0,1.5) -- (0,-1.5);
  \def\actionupper{(-1.5,1.5)
    .. controls (0.1,0.87571605) and (0.80,0.40) .. (0.80,-0.10)}
  \path[name path=action dotted] \actionupper;
  \path[name intersections={of=action dotted and action vertical,
                            by=action crossing}];
  \begin{scope}[even odd rule]
    \clip[overlay] (-2,-2) rectangle (2.5,2)
      (action crossing) circle (2.5pt);
    \draw[dotted stroke] \actionupper
      .. controls (0.80,-0.60) and (0.7,-1.05) .. (0,-1.05);
  \end{scope}
  \draw (0,1.5) -- (0,-1.5);
  \fill (0,-1.05) circle (0.4pt);
  \node[anchor=south, inner sep=2pt] at (0,1.5) {$\scriptstyle X$};
  \node[anchor=north, inner sep=2pt] at (0,-1.5) {$\scriptstyle X$};
\end{tikzpicture}%
}}\quad.
\]
These left and right $A$-action make $X$ to be an $A$-bimodule, denoted by $X^+$.

The induction functor $\a:\CC\to \CA$ is defined by $\a(x)=x\o A$. It is
both the left and the right adjoint to the forgetful functor $\CA\to \CC$
\cite[Lem. 1.16]{KO02}. Hence, for $x\in \CC$ and $Y\in \CA$, there is a
natural isomorphism $\CA(\a(x),Y)\cong \CC(x,Y)$. In particular, since $A$
is the unit object of $\CA$, we have $\CA(\a(x),\1) \stackrel{\vp}{\cong}\CC(x,A)$. In terms of diagrams,
\begin{equation}\label{eq:vp}
\vp_x(f)\,=\,
\vcenter{\hbox{%
\begin{tikzpicture}[x=1pc, y=1pc, thick]
  \node[draw, minimum width=20pt, minimum height=10pt,
        inner sep=0pt, outer sep=0pt] (f box) at (0,0) {$\scriptstyle f$};

  \node[inner sep=1pt, outer sep=0pt] (x input) at (-0.5,1.5) {$\scriptstyle x$};
  \draw (x input.south) -- (-0.5,0 |- f box.north);

  \node[inner sep=0pt, outer sep=0pt] (open input) at (0.5,1.5) {$\circ$};
  \draw[dotted stroke] (open input.south) -- (0.5,0 |- f box.north);

  \draw[dotted stroke] (f box.south) -- (0,-1.5);
\end{tikzpicture}%
}} \quad \text{ for } f \in \CA(\a(x), \1)\,.
\end{equation}
The inverse of $\vp_x$, denoted by $\psi_x$, is given by the diagram:
\begin{equation}\label{eq:psi}
   \psi_x(g)\,=\,
 \vcenter{\hbox{%
\begin{tikzpicture}[x=1pc, y=1pc, thick]
  \node[draw, minimum width=10pt, minimum height=10pt,
        inner sep=0pt, outer sep=0pt] (g box) at (0,0) {$\scriptstyle g$};

  \node[inner sep=1pt, outer sep=0pt] (x input) at (0,1.5) {$\scriptstyle x$};
  \draw (x input.south) -- (g box.north);

  \coordinate (dotted join) at (0,-1.5);
  \draw[dotted stroke] (g box.south) -- (0,-2);

  \draw[dotted stroke] (1,0 |- x input.south)
    .. controls (1,1) and (0.8,-1.2) .. (dotted join);
  \fill (dotted join) circle (0.4pt);
\end{tikzpicture}%
}} \quad \text{ for } g \in \CC(x, A)\,. 
\end{equation}

Let $\CAo$ be the category of local $A$-modules in $\CC$. More explicitly, for a
right $A$-module $Y$, the braiding of $\CC$ gives two left $A$-actions on
$Y$, and $Y$ is local if these two left actions coincide. The category
$\CAo$ is a fusion subcategory of $\CA$, and the braiding of $\CC$ descends
to a nondegenerate braiding on $\CAo$. Thus $\CAo$ is again an MTC
\cite[Thm. 4.5]{KO02}.

We will use the following result repeatedly. The functor $\a$ gives a fully
faithful braided tensor functor $\a^+:\CC\to Z(\CA)$. On the other hand,
local modules give a fully faithful braided tensor functor
$\b:(\CAo)^{\rev}\to Z(\CA)$. The two images centralize each other, and the
functor $\a^+\boxtimes \b:\CC\boxtimes(\CAo)^{\rev}\to Z(\CA)$, given by
$x\boxtimes X\mapsto \a^+(x)\o \b(X)$, is a braided tensor equivalence
preserving pivotal structures \cite[Cor. 3.28, Cor. 3.30]{DMNO}.

Since $\CC$ is semisimple, the algebra object $A$ has a decomposition
$$A=\bigoplus_{x\in\irr(\CC)} W_x\o x=\bigoplus_{x\in\irr(\CC)} V_x\o x,$$
where $W_x=\CC(x,A)$ is the
multiplicity space of $x$ in $A$ and $V_x=\CC_A(\alpha(x),A)$. The fusion action of $\CA$ on $W_x$ can be written explicitly as follows \cite{DNRX2025}.
For $x\in \CC$, $Y\in \CA$ and $g\in W_x=\CC(x,A)$, define
$Yg\in \CC(x,A)$ by
\begin{align*}
Yg=\Big( x &\xrightarrow{x\o \coev_Y} x\o Y\o Y^*
\xrightarrow{R_{x,Y}\o Y^*} Y\o x\o Y^*
\xrightarrow{Y\o g\o Y^*} Y\o A\o Y^*  \\
&\xrightarrow{\mu_Y\o Y^*} Y\o Y^*
=Y^{**}\o Y^*
\xrightarrow{\widetilde{\ev}_{Y^*}} A \Big),
\end{align*}
where $\mu_Y:Y\o A\to Y$ is the right $A$-module structure and
$\widetilde{\ev}_{Y^*}$ is the morphism descending to the right evaluation of $Y$ in $\CA$. 
In the graphical notation used below, this action is
\begin{equation}\label{eq:action_in_C}
Yg:=\frac{1}{d(A)}\vcenter{\hbox{%
\begin{tikzpicture}[x=1pc, y=1pc, thick]
  \draw ({1.5*cos(135)},{2.1*sin(135)})
    arc[start angle=135, end angle=476.56505,
        x radius=1.5pc, y radius=2.1pc];

  \node[draw, minimum width=10pt, minimum height=10pt,
        inner sep=0pt, outer sep=0pt] (Yg box) at (0,0.3)
    {$\scriptstyle g$};

  \draw (Yg box.north)
    .. controls (-0.1,1.5) and (-2,2) .. (-1.5,3.3);
  \node[anchor=south, inner sep=2pt] at (-1.5,3.3)
    {$\scriptstyle x$};

  \coordinate (Yg left join) at ({-1.5*sqrt(1-(1.1/2.1)^2)},-1.1);
  \coordinate (Yg right join) at ({1.5*sqrt(1-(1.3/2.1)^2)},-1.3);
  \draw[dotted stroke] (Yg box.south)
    .. controls (-0.2,-0.73333) and (-0.6,-1.03333)
    .. (Yg left join);
  \draw[dotted stroke] (Yg right join)
    .. controls (2.5,1) and (2.3,-1.7) .. (2.5,-2.8);
  \fill (Yg left join) circle (0.4pt);
  \fill (Yg right join) circle (0.4pt);
\end{tikzpicture}%
}}\,\text{, where}\qquad
\widetilde{\ev}_{Y^*}:=\frac{1}{d(A)}\vcenter{\hbox{%
\begin{tikzpicture}[x=1pc, y=1pc, thick]
  \draw (-1.5,0)
    arc[start angle=180, end angle=360,
        x radius=1.5pc, y radius=2.1pc];

  \node[anchor=south, inner sep=2pt] at (-1.5,0) {$\scriptstyle Y^{**}$};
  \node[anchor=south, inner sep=2pt] at (1.5,0) {$\scriptstyle Y^*$};

  \coordinate (Yg right join) at ({1.5*sqrt(1-(1.3/2.1)^2)},-1.3);
  \draw[dotted stroke] (Yg right join)
    .. controls (2.5,1) and (2.3,-1.7) .. (2.5,-2.8);
  \fill (Yg right join) circle (0.4pt);
\end{tikzpicture}%
}}\,.
\end{equation}
On the other hand, there is an action of $\CA$ on $V_x$ via the notion of generalized Frobenius-Schur indicators (cf. \cite[Sec. 2]{NS10}). The readers are referred to \cite{NSind, NS07, NS10, MNg, LM, KSZ, SZ, FGSV} for the development of Frobenius-Schur indicators.

    In terms of diagrams, the action of $Y \in \CA$ for $f\in V_x$ is given by
\begin{equation}\label{eq: indicator_diag}
Yf:=\vcenter{\hbox{%
\begin{tikzpicture}[x=1pc, y=1pc, line width=0.6pt]
  \draw[blue] ({1.3*cos(150.64649)},{1.5*sin(150.64649)})
    arc[start angle=150.64649, end angle=489.09356,
        x radius=1.3pc, y radius=1.5pc];

  \node[draw=blue, text=blue, minimum width=10pt, minimum height=10pt,
        inner sep=0pt, outer sep=0pt] (Yf box) at (0,0)
    {$\scriptstyle f$};

  \node[inner sep=1pt, outer sep=0pt] (Yf input) at (-1.5,2.2)
    {$\scriptstyle\alpha(x)$};
  \draw[blue] (Yf box.north)
    .. controls (-0.3,1) and (-2,0.7) .. (Yf input.south);

  \node[inner sep=1pt, outer sep=0pt] at (-1.8,-0.5)
    {$\scriptstyle Y$};
\end{tikzpicture}%
}}\hspace{0.5cm}.
\end{equation}
Here the blue lines mean the diagram is read in the category $\CA$, and the crossing in the diagram denotes the half-braiding of $\a^+(x)$.
By definition \cite[Sec. 2]{NS10} (cf. also \cite[Defn. 3.1]{DNRX2025}), $Y$ is the operator on $\CA(\a^+(x), \1)$, which defines the generalized Frobenius-Schur indicator $\nu_{(0,-1)}^{\a^+(x)}(Y)$.  More precisely, the ordinary trace of $Y$ on $V_x$ is given by
\begin{equation}\label{eq:GSFTrace}
    \Tr(Y, V_x) = \nu_{(0,-1)}^{\a^+(x)}(Y) \,.
\end{equation}
It follows from \cite{DNRX2025} that these actions make $W_x$ and $V_x$ being equivalent $K(\CA)$-modules, and so we find
\begin{equation}\label{eq:GSFTrace}
    \Tr(Y, W_x) = \nu_{(0,-1)}^{\a^+(x)}(Y) \,.
\end{equation}

In particular, by taking $x=A$ and $g=\id_A$ in \eqref{eq:action_in_C}, we obtain the fusion action
of $\CA$ on the algebra  $A$. Equivalently, it is the $\BC$-linear
map $\rho:K(\CA)\to \CC(A,A)$ defined on simple objects by
\begin{equation}\label{eq:rho}
\rho(Y)=Y\id_A=\frac{1}{d(A)}\vcenter{\hbox{%
\begin{tikzpicture}[x=0.9pc, y=0.9pc, thick]
  \draw ({1.5*cos(66.36419)},{2.2*sin(66.36419)})
    arc[start angle=66.36419, end angle=413.12867,
        x radius=1.35pc, y radius=1.98pc];

  \coordinate (rho left join) at ({-1.5*sqrt(1-(0.3/2.2)^2)},-0.3);
  \draw[dotted stroke] (1.3,3)
    .. controls (1.03333,2) and (0.10467,0.9)
    .. (rho left join);

  \coordinate (rho right join) at ({1.5*sqrt(1-(1.3/2.2)^2)},-1.3);
  \draw[dotted stroke] (rho right join)
    .. controls (2.5,1) and (2.3,-1.7) .. (2.4,-2.8);
  \fill (rho left join) circle (0.4pt);
  \fill (rho right join) circle (0.4pt);

  \node[inner sep=1pt, outer sep=0pt] at (-1.4,1.5)
    {$\scriptstyle Y$};
\end{tikzpicture}%
}}\quad\text{for }Y\in\iCA.
\end{equation}
The associativity of the action defined in \eqref{eq:action_in_C} implies that $\rho: K(\CA)\to \CC(A,A)$ is an algebra
homomorphism, and this homomorphism $\rho$ will be called 
the fusion action of $\CA$ on $A$.

For any fusion subcategory $\BB$ of $\CA$, we will write
$\CC(x,A)^\BB$ for the $\BB$-invariant subspace of $\CC(x,A)$, namely
$$
\CC(x,A)^\BB=\{g\in \CC(x,A)\mid Yg=d_A(Y)g \text{ for all }Y\in\BB\}.
$$
This notation will be used to define fixed-point  of $A$ under
fusion actions.

\section{Schur--Weyl duality and Galois correspondence}

In this section, we recall two results from our previous paper
\cite{DNRX2025} with some improvement on the second one. They will be used as the main categorical input in the present
paper.

Let $\CC$ be an MTC and let $A$ be a condensable algebra in $\CC$. Let
$\CA$ denote the fusion category of right $A$-modules in $\CC$, and  $\CAo$
 the fusion subcategory of local $A$-modules. Recall that $\CAo$ is an
MTC. The fusion action constructed in \cite{DNRX2025} gives a $\BC$-algebra
homomorphism
$\rho:K(\CA)\to \CC(A,A)$. For $Y\in\irr(\CA)$, $\rho(Y)$ is the
endomorphism of $A$ obtained by the action of $Y$ on
$\id_A\in\CC(A,A)$.

We shall use the following idempotent associated with the dimension
character of $K(\CAo)$:
\[
e_\1=\frac{1}{\dim(\CAo)}\sum_{X\in\irr(\CAo)}d_A(X)X.
\]
Since $K(\CAo)$ is a central subalgebra of $K(\CA)$,  $e_\1$ as a central idempotent
of $K(\CA)$.

The first result is the Schur--Weyl duality for the fusion action of
$\CA$ on $A$.

\begin{thm}[Schur--Weyl duality, {\cite[Theorem I]{DNRX2025}}]\label{t:SW-fusion}
Let $\CC$ be an MTC and let $A$ be a condensable algebra in $\CC$. Let
$\rho:K(\CA)\to\CC(A,A)$ be the fusion action. Then:
\begin{enumerate}
\item $\ker\rho=(1-e_\1)K(\CA)$.
\item For any $x\in\irr(\CC)$ with $\CC(x,A)\ne 0$, the space $\CC(x,A)$ is
an irreducible $K(\CA)$-module.
\item For $x,y\in\irr(\CC)$ with $\CC(x,A)\ne 0$ and $\CC(y,A)\ne 0$, one
has $\CC(x,A)\cong\CC(y,A)$ as $K(\CA)$-modules if and only if $x=y$.
\item The restriction $\rho:e_\1K(\CA)\to\CC(A,A)$ is an isomorphism of
$\BC$-algebras.
\end{enumerate}
\end{thm}

Equivalently, if $A=\bigoplus_{x\in\irr(\CC)}W_x\o x$ with
$W_x=\CC(x,A)$, then the nonzero multiplicity spaces $W_x$ are precisely
the irreducible modules for the effective fusion algebra $e_\1K(\CA)$
which occur in $A$. Thus the fusion action gives a categorical analogue of
the Schur--Weyl decomposition in orbifold theory.

We next recall the fixed-point construction for fusion subcategories.   Let
$\BB$ be a fusion subcategory of $\CA$ containing $\CAo$. Put
\[
e_\1^\BB=\frac{1}{\dim(\BB)}\sum_{Y\in\irr(\BB)}d_A(Y)Y.
\]
Then $\rho(e_\1^\BB)$ is an idempotent in $\CC(A,A)$. Let
$\rho(e_\1^\BB)=\iota\pi$ be its epi-monic factorization, where
$\iota:B\to A$ is the image and $\pi:A\to B$ is the coimage. We denote this
subobject of $A$ by $A^\BB=(B,\iota)$. The subobject $B$ inherits an algebra structure so that $\iota: B \to A$ is an algebra monomorphism or $(B,\iota)$ is a subalgebra of $A$, i.e. $\iota m_B =m_A (\iota \o  \iota)$ and $\iota u_B = u_A$. We will call $(B,\iota)$ the \emph{fixed points} of $A$ under the action of $\BB.$ In fact, as an object of $\CC$,
\[
A^\BB=\sum_{x\in\irr(\CC)}\CC(x,A)^\BB\o x.
\]

Secondly, we  strengthen the Galois correspondence result established in \cite[Theorem II]{DNRX2025} without the pseudounitary assumption.

\begin{lem}\label{l:con_sub}
 Every subalgebra $B$ of $A$ with $d(B) \ne 0$ is a   condensable algebra in $\CC$.
 \end{lem} 
 \begin{proof}
 Since  $B$ is a subalgebra of $A$, $B$ is connected and commutative with $\theta_B = \id_B$. Under the assumption that $d(B)\ne 0$, it suffices to show that the algebra $B$ is simple as $B$-module by \cite[Lem. 1.20]{KO02}. However,  the isomorphisms \cite[Thm. 1.6]{KO02}
 $$
\CB(B,B) \cong \CC(\1, B) \cong \BC
 $$  
 implies $B$ is a simple $B$-module, and hence $B$ is a condensable algebra in $\CC$. 
 \end{proof}
 \begin{prop} \label{p:positivity} Let $A$ be a condensable algebra in a modular fusion category $\CC$. Then
 \begin{enumeri}
     \item Every simple subobject $x$ of  $A$  has totally positive categorical dimension, which means $\sigma(d(x)) > 0$ for all automorphisms $\sigma  \in \GQ$.
     \item Every subalgebra $B$ of $A$ is condensable and has  totally positive dimension.
 \end{enumeri}
 \end{prop}
 \begin{proof}
     (i) Recall that $Z(\CA) \simeq \CC \boxtimes (\CAo)^{\rev}$ as modular tensor categories and $ A \cong A \boxtimes\1$ is a subalgebra of the condensable algebra $I(\1) \in Z(\CA)$ as $A = I(\1)^{\CAo}$ (cf. \cite[Prop. 5.5]{DNRX2025}). It suffices to show that each simple subobject of $I(\1)$ has totally positive categorical dimension. We may simply assume $A = I(\1)$.

     By Theorem \ref{t:SW-fusion}, for any simple subobject $x$ of $A$, $W_x= \CC(x,A) \cong \CA(\a(x), \1)$ is an  irreducible representation of $K(\CA)$ with the character
     $\chi_{W_x}$ given by
      the generalized FS-indicator (cf. \eqref{eq:GSFTrace})
      $$
      \chi_{W_x}(Y) = \nu_{(0,-1)}^{\a^+(x)}(Y)\,.
      $$
      So, by \cite[Lem. 5.3 and Cor. 5.6]{NS10} we find
      $$
      \chi_{W_x}(Y) = \nu_{(0,1)}^{\a^+(x)^*}(Y) = \frac{1}{\dim(\CA)} \sum_{z \in \irr(\CC)} S_{x, z}[Y,x]_\CC\,.
      $$
      It follows from \cite[Prop. 3.3]{NS10} that
      \begin{align*}
          \chi_{W_x}(Y^*) & = \nu_{(0,1)}^{\a^+(x)}(Y)=\frac{1}{\dim(\CA)} \sum_{z \in \irr(\CC)} S_{x^*, z}[Y,x]_\CC \\
          &=\frac{1}{\dim(\CA)} \sum_{z \in \irr(\CC)} \ol{S}_{x, z}[Y,x]_\CC 
          =\ol{\chi_{W_x}(Y)} \,.
      \end{align*}
     Thus, the formal codegree  $f_{W_x}$ of $W_x$ given by  
     $$
      \frac{1}{\dim_\BC W_x} \sum_{Y \in \iCA} \chi_{W_x}(Y) \chi_{W_x}(Y^*) =\frac{1}{\dim_\BC W_x} \sum_{Y \in \iCA} |\chi_{W_x}(Y)|^2
     $$   
     is totally positive. By \cite[Thm. 2.13]{O15} or the proof of \cite[Thm. 4.1 (ii)]{DNRX2025} that
     $$
     f_{W_x} = \frac{\dim(\CA)}{d_A(\a(x))} =  \frac{\dim(\CA)}{d(x)}\,.
     $$
     Since $\dim(\CA)$ is also totally positive, and so is $d(x)$.
     
(ii) follows directly from Lemma \ref{l:con_sub} and (i).
 \end{proof}

 In part (iv) of the proof of \cite[Theorem II]{DNRX2025}, the pseudounitarity of $\CC$ was used to prove the following lemma. Here, we provide a proof without such assumption.
 \begin{lem} \label{l:3.4}
 Let $\CC$ be a modular fusion category, $A$ a condensable algebra in $\CC$, $\BB$ a fusion category of $\CA$ containing $\CAo$, and $(B, \iota)=A^\BB$. Then every object $X \in \BB$ considered as a $B$-module via the embedding $\iota$ is a local $B$-module, or $\BB \subseteq (\CBo)_A$. 
 \end{lem}
 \begin{proof} Let $X \in \irr(\BB)$. In part (iv) of the proof of \cite[Theorem II]{DNRX2025}, without using pseudounitarity of $\CC$, it has established that $\tr_\CC(P_{X}) = d(X)$ where 
 \begin{equation}
P_{X} = \frac{1}{d(B)}\begin{tikzpicture}[scale=1,
  x=0.2cm, y=0.2cm,
  thick,
  baseline={([yshift=-\the\fontdimen22\textfont2]diagram midpoint)}
]
  \coordinate[overlay] (diagram midpoint) at (0,5.375);
  \draw (0,10) -- (0,5.45);
  \draw (0,4.30) -- (0,-1);

  \draw
    (-0.500,8.020)
    .. controls (-1.167,7.920) and (-1.872,6.520) .. (-1.550,5.700)

    .. controls (-1.228,4.880) and (0.123,4.639) .. (0.900,4.300)

    .. controls (1.677,3.961) and (1.882,3.525) .. (1.900,2.950)

    .. controls (1.918,2.75) and (1.750,1.660) .. (0.000,1.660);

  \draw
    (0.50,8.15)
    .. controls (2.50,8.15) and (4.10,6.40) .. (4.10,4.245)
    .. controls (4.10,2.09) and (2.20,0.34) .. (0.00,0.34);

  \node[anchor=south, inner sep=2pt] at (0,10) {$\scriptstyle X$};
  \node[anchor=west, inner sep=3pt] at (4.10,4.245) {$\scriptstyle B$};
\end{tikzpicture}
 \end{equation}
 is a projection onto the largest local $B$-submodule  $X'$ of $X$. It suffices to show that $P_X$ is an $A$-module morphism without assuming pseudounitarity of $\CC$, and so $P_X = \lambda \id_X$. The equality $\tr_\CC(P_{X}) = d(X)$ implies $P_X = \id_X$ or $X = X'$. 

 We proceed to prove $P_X$ is an $A$-module morphism by graphical calculus. Since $A$ is commutative, and $B$ acts on $X$ and $A$ via $\iota$, we have
\[\displaystyle
\begin{tikzpicture}[x=1cm, y=1cm, xscale=0.75, yscale=0.6, thick, baseline={([yshift=-0.5ex]diagram midpoint)}]
  \coordinate[overlay] (diagram midpoint) at (0,1.90098);
  \draw (0,0) -- (0,3.80196);

  \draw[name path=solid strand] (0,0.97490)
    .. controls (1.03959,1.53307) and (1.01337,2.45823)
    .. (0.92134,3.75036);

  \def\dottedcurve{(0,1.68620)
    .. controls (1.17127,1.90698) and (1.42270,2.59039)
    .. (1.75428,3.73635)}
  \path[name path=dotted strand] \dottedcurve;

  \path[name intersections={of=solid strand and dotted strand,by=crossing}];
  \begin{scope}[even odd rule]
    \clip[overlay] (-0.5,-0.5) rectangle (3,4.5)
      (crossing) circle (1.5mm);
    \draw[dotted stroke]
      \dottedcurve;
  \end{scope}

  \node[anchor=south, inner sep=3pt ] at (0,3.80196) {$\scriptstyle X$};
  \node[anchor=south, inner sep=3pt ] at (0.92134,3.75036) {$\scriptstyle B$};
\end{tikzpicture}
=
\begin{tikzpicture}[x=1cm, y=1cm, xscale=0.75, yscale=0.6, thick, baseline={([yshift=-0.5ex]diagram midpoint)}]
  \coordinate[overlay] (diagram midpoint) at (0,1.90098);
  \draw (0,0) -- (0,3.80196);

  \def\bcurve{(0,0.97490)
    .. controls (1.03959,1.53307) and (1.01337,2.45823)
    .. (0.92134,3.75036)}
  \path[name path=B strand] \bcurve;

  \def\dottedcurve{(0,1.68620)
    .. controls (1.17127,1.90698) and (1.42270,2.59039)
    .. (1.75428,3.73635)}
  \path[name path=dotted strand] \dottedcurve;

  \path[name intersections={of=B strand and dotted strand,by=crossing}];
  \begin{scope}[even odd rule]
    \clip[overlay] (-0.5,-0.5) rectangle (3,4.5)
      (crossing) circle (1.5mm);
    \draw[dotted stroke] \dottedcurve;
  \end{scope}

  \path[overlay, name path=iota level] (0,3.25036) -- (2,3.25036);
  \path[name intersections={of=B strand and iota level,by=iota position}];

  \begin{scope}
    \clip[overlay] (-0.5,-0.5) rectangle (3,0 |- iota position);
    \draw[dotted stroke] \bcurve;
  \end{scope}
  \begin{scope}
    \clip[overlay] (-0.5,0 |- iota position) rectangle (3,4.5);
    \draw \bcurve;
  \end{scope}

  \node[draw, fill=white, line width=0.6pt, minimum size=4mm, inner sep=0pt]
    at (iota position) {$\scriptstyle\iota$};

  \node[anchor=south, inner sep=3pt ] at (0,3.80196) {$\scriptstyle X$};
  \node[anchor=south, inner sep=3pt ] at (0.92134,3.75036) {$\scriptstyle B$};
\end{tikzpicture}
=
\begin{tikzpicture}[x=1cm, y=1cm, xscale=0.75, yscale=0.6, thick, baseline={([yshift=-0.5ex]diagram midpoint)}]
  \coordinate[overlay] (diagram midpoint) at (0,1.90098);
  \draw (0,0) -- (0,3.80196);

  \def\bcurve{(0,0.97490)
    .. controls (1.03959,1.53307) and (1.01337,2.45823)
    .. (0.92134,3.75036)}
  \path[name path=B strand] \bcurve
    coordinate[pos=0.1324788855] (B join);

  \coordinate (dotted turn) at (0.65238112,1.87418296);
  \def\dotteduppercontinuation{
    .. controls (1.455048,2.198846) and (1.489016,2.819582)
    .. (1.75428,3.73635)}

  \def\dottedcurve{(B join)
    .. controls (0.42000,1.55000) and (0.46000,1.79637)
    .. (dotted turn)
    \dotteduppercontinuation}

  \path[name path=dotted strand] (dotted turn) \dotteduppercontinuation;

  \path[name intersections={of=B strand and dotted strand,by=crossing}];
  \begin{scope}[even odd rule]
    \clip[overlay] (-0.5,-0.5) rectangle (3,4.5)
      (crossing) circle (1.5mm);
    \draw[dotted stroke] \dottedcurve;
  \end{scope}

  \path[overlay, name path=iota level] (0,3.25036) -- (2,3.25036);
  \path[name intersections={of=B strand and iota level,by=iota position}];

  \begin{scope}
    \clip[overlay] (-0.5,-0.5) rectangle (3,0 |- iota position);
    \draw[dotted stroke] \bcurve;
  \end{scope}
  \begin{scope}
    \clip[overlay] (-0.5,0 |- iota position) rectangle (3,4.5);
    \draw \bcurve;
  \end{scope}

  \fill (B join) circle (0.4pt);

  \node[draw, fill=white, line width=0.6pt, minimum size=4mm, inner sep=0pt]
    at (iota position) {$\scriptstyle\iota$};

  \node[anchor=south, inner sep=3pt ] at (0,3.80196) {$\scriptstyle X$};
  \node[anchor=south, inner sep=3pt ] at (0.92134,3.75036) {$\scriptstyle B$};
\end{tikzpicture}
=
\begin{tikzpicture}[x=1cm, y=1cm, xscale=0.75, yscale=0.6, thick, baseline={([yshift=-0.5ex]diagram midpoint)}]
  \coordinate[overlay] (diagram midpoint) at (0,1.90098);
  \draw (0,0) -- (0,3.80196);

  \def\bcurve{(0,0.97490)
    .. controls (1.03959,1.53307) and (1.01337,2.45823)
    .. (0.92134,3.75036)}
  \path[name path=B strand] \bcurve
    coordinate[pos=0.35] (B join);

  \draw[dotted stroke] (B join)
    .. controls ([xshift=3mm,yshift=2mm]B join) and (2.10000,2.75000)
    .. (1.75428,3.73635);

  \path[overlay, name path=iota level] (0,{3.25036-1/3}) -- (2,{3.25036-1/3});
  \path[name intersections={of=B strand and iota level,by=iota position}];

  \begin{scope}
    \clip[overlay] (-0.5,-0.5) rectangle (3,0 |- iota position);
    \draw[dotted stroke] \bcurve;
  \end{scope}
  \begin{scope}
    \clip[overlay] (-0.5,0 |- iota position) rectangle (3,4.5);
    \draw \bcurve;
  \end{scope}

  \fill (B join) circle (0.4pt);

  \node[draw, fill=white, line width=0.6pt, minimum size=4mm, inner sep=0pt]
    at (iota position) {$\scriptstyle\iota$};

  \node[anchor=south, inner sep=3pt ] at (0,3.80196) {$\scriptstyle X$};
  \node[anchor=south, inner sep=3pt ] at (0.92134,3.75036) {$\scriptstyle B$};
\end{tikzpicture}
=
\begin{tikzpicture}[x=1cm, y=1cm, xscale=0.75, yscale=0.6, thick, baseline={([yshift=-0.5ex]diagram midpoint)}]
  \coordinate[overlay] (diagram midpoint) at (0,1.90098);
  \draw (0,0) -- (0,3.80196);

  \def\bcurve{(0,0.97490)
    .. controls (1.03959,1.53307) and (1.01337,2.45823)
    .. (0.92134,3.75036)}
  \path[name path=B strand] \bcurve;

  \draw[dotted stroke] (0,0.47490)
    .. controls (1.15000,0.95000) and (2.10000,2.75000)
    .. (1.75428,3.73635);

  \path[overlay, name path=iota level] (0,{3.25036-1/3}) -- (2,{3.25036-1/3});
  \path[name intersections={of=B strand and iota level,by=iota position}];

  \begin{scope}
    \clip[overlay] (-0.5,-0.5) rectangle (3,0 |- iota position);
    \draw[dotted stroke] \bcurve;
  \end{scope}
  \begin{scope}
    \clip[overlay] (-0.5,0 |- iota position) rectangle (3,4.5);
    \draw \bcurve;
  \end{scope}

  \node[draw, fill=white, line width=0.6pt, minimum size=4mm, inner sep=0pt]
    at (iota position) {$\scriptstyle\iota$};

  \node[anchor=south, inner sep=3pt ] at (0,3.80196) {$\scriptstyle X$};
  \node[anchor=south, inner sep=3pt ] at (0.92134,3.75036) {$\scriptstyle B$};
\end{tikzpicture}
= 
\begin{tikzpicture}[x=1cm, y=1cm, xscale=0.75, yscale=0.6, thick, baseline={([yshift=-0.5ex]diagram midpoint)}]
  \coordinate[overlay] (diagram midpoint) at (0,1.90098);
  \draw (0,0) -- (0,3.80196);

  \draw (0,0.97490)
    .. controls (1.03959,1.53307) and (1.01337,2.45823)
    .. (0.92134,3.75036);

  \draw[dotted stroke] (0,0.47490)
    .. controls (1.15000,0.95000) and (1.82453,2.75000)
    .. (1.75428,3.73635);

  \node[anchor=south, inner sep=3pt ] at (0,3.80196) {$\scriptstyle X$};
  \node[anchor=south, inner sep=3pt ] at (0.92134,3.75036) {$\scriptstyle B$};
\end{tikzpicture}\,.
\]
The diagram equalities obviously hold for the reversed braids, and so 
 \[
\begin{tikzpicture}[scale=1,
  x=0.2cm, y=0.2cm,
  thick,
  baseline={([yshift=-\the\fontdimen22\textfont2]diagram midpoint)}
]
  \coordinate[overlay] (diagram midpoint) at (0,5.375);
  \draw (0,11.75) -- (0,5.45);
  \draw (0,4.30) -- (0,-1);

  \draw[dotted stroke] (0,9.02)
    .. controls (1.40,9.78) and (3.40,11.28)
    .. (3.60,12.8);

  \draw
    (-0.500,8.020)
    .. controls (-1.167,7.920) and (-1.872,6.520) .. (-1.550,5.700)

    .. controls (-1.228,4.880) and (0.123,4.639) .. (0.900,4.300)

    .. controls (1.677,3.961) and (1.882,3.525) .. (1.900,2.950)

    .. controls (1.918,2.75) and (1.750,1.660) .. (0.000,1.660);

  \draw
    (0.50,8.15)
    .. controls (2.50,8.15) and (4.10,6.40) .. (4.10,4.245)
    .. controls (4.10,2.09) and (2.20,0.34) .. (0.00,0.34);

  \node[anchor=south, inner sep=2pt] at (0,11.75) {$\scriptstyle X$};
  \node[anchor=west, inner sep=3pt] at (4.10,4.245) {$\scriptstyle B$};
\end{tikzpicture}
\;=\;
\begin{tikzpicture}[scale=1,
  x=0.2cm, y=0.2cm,
  thick,
  baseline={([yshift=-\the\fontdimen22\textfont2]diagram midpoint)}
]
  \coordinate[overlay] (diagram midpoint) at (0,4.375);
  \draw (0,10.75) -- (0,5.45);
  \draw (0,4.30) -- (0,-2);

  \node[anchor=south, inner sep=2pt] (x label) at (0,10.75) {$\scriptstyle X$};
  \coordinate (dotted top) at (3.60,0 |- x label.center);

  \def\dottedcurve{(0,3)
    .. controls (2.5,4) and (3.40,11.28)
    .. (dotted top)}
  \path[name path=dotted strand] \dottedcurve;

  \draw[name path=left strand]
    (-0.500,8.020)
    .. controls (-1.167,7.920) and (-1.872,6.520) .. (-1.550,5.700)

    .. controls (-1.228,4.880) and (0.123,4.639) .. (0.900,4.300)

    .. controls (1.677,3.961) and (1.882,3.525) .. (1.900,2.950)

    .. controls (1.918,2.75) and (1.750,1.660) .. (0.000,1.660);

  \path[name intersections={of=dotted strand and left strand,by=left crossing}];
  \begin{scope}[even odd rule]
    \clip
      (-3,-1) rectangle (6,14)
      (left crossing) circle (1.6mm);
    \draw[dotted stroke] \dottedcurve;
  \end{scope}

  \path[name path=right strand]
    (0.50,8.15)
    .. controls (2.50,8.15) and (4.10,6.40) .. (4.10,4.245)
    .. controls (4.10,2.09) and (2.20,0.34) .. (0.00,0.34);

  \path[name intersections={of=dotted strand and right strand,by=right crossing}];
  \begin{scope}[even odd rule]
    \clip
      (-3,-1) rectangle (6,14)
      (right crossing) circle (2.2mm);
    \draw
      (0.50,8.15)
      .. controls (2.50,8.15) and (4.10,6.40) .. (4.10,4.245)
      .. controls (4.10,2.09) and (2.20,0.34) .. (0.00,0.34);
  \end{scope}

  \node[anchor=west, inner sep=3pt] at (4.10,4.245) {$\scriptstyle B$};
\end{tikzpicture}
\;=\;
\begin{tikzpicture}[
  x=0.2cm, y=0.2cm,
  thick,
  baseline={([yshift=-\the\fontdimen22\textfont2]diagram midpoint)}
]
  \coordinate[overlay] (diagram midpoint) at (0,4.375);
  \path[name path=vertical strand] (0,10.75) -- (0,-2);

  \node[anchor=south, inner sep=2pt] (x label) at (0,10.75) {$\scriptstyle X$};
  \coordinate (dotted top) at (3.60,0 |- x label.center);

  \def\dottedcurve{(0,1)
    .. controls (4.0,4) and (3.40,11.28)
    .. (dotted top)}
  \path[name path=dotted strand] \dottedcurve;

  \def\leftuppercurve{
    (-0.500,8.020)
    .. controls (-1.167,7.920) and (-1.872,6.520) .. (-1.550,5.700)

    .. controls (-1.228,4.880) and (0.123,4.639) .. (0.900,4.300)

    .. controls (1.100,4.210) and (1.100,4.000) .. (1.100,3.750)}
  \path[name path=left strand] \leftuppercurve;

  \path[name intersections={of=vertical strand and left strand,by=vertical crossing}];
  \draw (0,10.75) -- ([yshift=1.15mm]vertical crossing);
  \draw ([yshift=-1.15mm]vertical crossing) -- (0,-2);

  \draw \leftuppercurve
    .. controls (1.100,3.050) and (0.650,2.400) .. (0,2.400);
  \draw[dotted stroke] \dottedcurve;

  \path[name path=right strand]
    (0.50,8.15)
    .. controls (2.50,8.15) and (4.10,6.40) .. (4.10,4.245)
    .. controls (4.10,2.09) and (2.20,-0.16) .. (0.00,-0.16);

  \path[name intersections={of=dotted strand and right strand,by=right crossing}];
  \begin{scope}[even odd rule]
    \clip
      (-3,-1) rectangle (6,14)
      (right crossing) circle (2.2mm);
    \draw
      (0.50,8.15)
      .. controls (2.50,8.15) and (4.10,6.40) .. (4.10,4.245)
      .. controls (4.10,2.09) and (2.20,-0.16) .. (0.00,-0.16);
  \end{scope}

  \node[anchor=west, inner sep=3pt] at (4.10,4.245) {$\scriptstyle B$};
\end{tikzpicture}
\;=\;
\begin{tikzpicture}[
  x=0.2cm, y=0.2cm,
  thick,
  baseline={([yshift=-\the\fontdimen22\textfont2]diagram midpoint)}
]
  \coordinate[overlay] (diagram midpoint) at (0,3.375);
  \path[name path=vertical strand] (0,9.75) -- (0,-3.00);

  \node[anchor=south, inner sep=2pt] (x label) at (0,9.75) {$\scriptstyle X$};
  \coordinate (dotted top) at (4.10,0 |- x label.center);

  \def\dottedcurve{(0,-2)
    .. controls (6.20,1.70) and (4.10,10.10)
    .. (dotted top)}

  \def\leftuppercurve{
    (-0.500,8.020)
    .. controls (-1.167,7.920) and (-1.872,6.520) .. (-1.550,5.700)

    .. controls (-1.228,4.880) and (0.123,4.639) .. (0.900,4.300)

    .. controls (1.100,4.210) and (1.100,4.000) .. (1.100,3.750)}
  \path[name path=left strand] \leftuppercurve;

  \path[name intersections={of=vertical strand and left strand,by=vertical crossing}];
  \draw (0,9.75) -- ([yshift=1.15mm]vertical crossing);
  \draw ([yshift=-1.15mm]vertical crossing) -- (0,-3.00);

  \draw \leftuppercurve
    .. controls (1.100,3.050) and (0.650,2.400) .. (0,2.400);
  \draw[dotted stroke] \dottedcurve;

  \coordinate (right bend) at (2.750,4.300);
  \draw
    (0.50,8.15)
    .. controls (1.900,8.150) and (2.750,6.400) .. (right bend)
    .. controls (2.750,2.100) and (1.500,0.500) .. (0,0.500);

  \node[anchor=east, inner sep=2pt] at (right bend) {$\scriptstyle B$};
\end{tikzpicture}
\,. 
 \]
 Therefore, $P_X$ is an $A$-module morphism.
 \end{proof}
Now we can improve Galois correspondence result \cite[Theorem II]{DNRX2025} using Proposition \ref{p:positivity} and Lemma \ref{l:3.4}.
 
\begin{thm}[Galois correspondence, {\cite[Theorem II]{DNRX2025}}]\label{t:Galois-fusion}
Let $\CC$ be any MTC, and let $A$ be a condensable algebra in
$\CC$. Then the assignment
$\BB\mapsto A^\BB$ defines an order-reversing bijection between fusion
subcategories $\BB$ of $\CA$ containing $\CAo$ and condensable subalgebras
of $A$. The inverse assignment is given by
$(B,\iota)\mapsto(\CC_B^0)_A$. In particular, for every fusion subcategory
$\BB$ of $\CA$ containing $\CAo$, one has
\[
\BB=(\CC_{A^\BB}^0)_A,
\qquad
\dim(\BB)=\frac{\dim(\CC)}{d(A)d(A^\BB)}.
\]
\end{thm}

Thus the assignment $\BB\mapsto A^\BB$ is a categorical fixed-point
construction. A larger fusion subcategory of $\CA$ gives a smaller
condensable subalgebra of $A$, as in the classical Galois correspondence
and in orbifold theory. 

In the special case where $\CA$ is a faithful $G$-extension of $\CAo$, the
fusion action reduces to an ordinary finite group action. Namely, if
$\CA=\oplus_{g\in G}\CA(g)$ with $\CA(1)=\CAo$, then
\[
\ol g=\frac{1}{\dim(\CAo)}
\sum_{X\in\irr(\CA(g))}d_A(X)X
\]
embeds $G$ into $K(\CA)$, and $e_\1K(\CA)$ is naturally identified with
$\BC[G]$. Under this identification, the fusion action on $A$ becomes the
corresponding $G$-action on $A$.

\section{Automorphism groups of condensable algebras}

In this section, we study algebra automorphisms of any condensable algebra $A$ in a modular tensor category $\CC$. We prove that $\Auta(A)$ is finite and investigate the decomposition of $A$ and its fixed-point subalgebra under a subgroup of algebra automorphisms.

\begin{defn} Let $A$ be a condensable algebra in an MTC $\CC$.
    An isomorphism $\s \in \CC(A,A)$ is called an algebra automorphism of $A$ if
    $\s u_A  = u_A$ and $m_A (\s \o \s) = \s m_A$. The set of all algebra automorphisms of $A$ is denoted by $\Auta(A)$, and reserve the notation $\AutC(A)$ to denote the group of  automorphisms on $A$ in $\CC$. Obviously, $\Auta$ is a subgroup of $\AutC(A)$ with the product given by the composition of morphisms. 
\end{defn}
\begin{rmk}
An algebra endomorphism of $A$ preserves $u_A$ and its  kernel is a right $A$-submodules of $A$. Since $A$ is a simple $A$-module,  algebra endomorphisms of $A$ are isomorphisms.  
\end{rmk}
We will prove that $\Auta(A)$ is a finite group. To this end, we first establish the following lemma, which was proved in the VOA setting \cite{McRae2021} under the assumption $\dim A=|G|$; however, the argument applies in a more general setting.
\begin{lem}\label{p5.26} Let $\CC$ be a ribbon finite tensor category (not necessarily semisimple), and $A$ a commutative algebra in $\CC$ such that $A$ is simple (as an $A$-module), $\CC(\1, A)$ is 1-dimensional, and $\1$ is a direct summand of $A$. 
Then:
\begin{enumerate}
\item For any $X \in \CA$ and $g \in \End_\CC(X)$, $E(g) : X \to X$ defined by the diagram 
$$
E(g)=\vcenter{\hbox{\begin{tikzpicture}[squarednode/.style={rectangle, draw=black,  minimum size=1pt}, baseline=0ex]
\node[squarednode](A) at (0,0) { $\hspace{2pt}\scriptstyle g \hspace{2pt}$}; 
\draw[thick] (A.90)-- +(0,1) node [xshift=5pt] {{\tiny $X$}};
\draw[thick] (A.270)-- +(0,-1)  node [xshift=3pt]{{\tiny $X$}};
\draw[thick,densely dotted] (0,.5) .. controls +(.7,.8) and +(1,.2) .. (0,-.7);
\end{tikzpicture}}}
$$
is an $A$-module morphism.
\item For any $g \in \Auta(A)$, $\tr_\CC(g) = \delta_{g, 1} d(A)$. 
\end{enumerate} 
\end{lem}
\begin{proof}
Let $\e_A : A \to \1$ such that $\e_A u_A = \id_\1$. The simplicity of $A$ implies the pairing $\e_A m_A: A \o A \to \1$ is nondegenerate, and so one can choose $\ev_A = \e_A m_A$ (cf. \cite{KO02}). The associativity of $m_A$ implies that $A$ is a Frobenius algebra. Thus
\[
\vcenter{\hbox{%
\begin{tikzpicture}[
  x=1cm, y=1cm, thick, baseline=0ex,
  squarednode/.style={rectangle, draw=black, minimum size=1pt}
]
  \node[squarednode] (A) at (0,0)
    {$\scriptstyle\hspace{2pt}g\hspace{2pt}$};

  \draw (A.north) -- +(0,1)
    node[xshift=5pt] {\tiny $X$};
  \draw (A.south) -- +(0,-1)
    node[xshift=3pt] {\tiny $X$};

  \draw[dotted stroke] (0,0.5)
    .. controls +(0.7,0.8) and +(1,0.2) .. (0,-0.7);

  \draw[dotted stroke] (0.5,1.2)
    .. controls (0.35,1.00) and (0.17,0.82) .. (0,0.7);
\end{tikzpicture}%
}}
\quad=\quad
\vcenter{\hbox{%
\begin{tikzpicture}[
  x=1cm, y=1cm, thick, baseline=0ex,
  squarednode/.style={rectangle, draw=black, minimum size=1pt}
]
  \node[squarednode] (A) at (0,0)
    {$\scriptstyle\hspace{2pt}g\hspace{2pt}$};

  \draw (A.north) -- +(0,1)
    node[xshift=5pt] {\tiny $X$};
  \draw (A.south) -- +(0,-1)
    node[xshift=3pt] {\tiny $X$};

  \draw[dotted stroke] (0,0.5)
    .. controls +(0.7,0.8) and +(1,0.2) .. (0,-0.7)
    coordinate[pos=0.073569435] (dotted join);

  \draw[dotted stroke] (0.5,1.2)
    .. controls (0.399214,0.978674) and (0.268428,0.777348)
    .. (dotted join);
  \fill (dotted join) circle (0.4pt);
\end{tikzpicture}%
}}
\quad=\quad
\vcenter{\hbox{%
\begin{tikzpicture}[
  x=1cm, y=1cm, thick, baseline=0ex,
  squarednode/.style={rectangle, draw=black, minimum size=1pt}
]
  \node[squarednode] (A) at (0,0)
    {$\scriptstyle\hspace{2pt}g\hspace{2pt}$};

  \draw (A.north) -- +(0,1)
    node[xshift=5pt] {\tiny $X$};
  \draw (A.south) -- +(0,-1)
    node[xshift=3pt] {\tiny $X$};

  \draw[dotted stroke] (0,0.5)
    .. controls +(0.7,0.8) and +(1,0.2) .. (0,-0.7)
    coordinate[pos=0.45] (dotted join);

  \draw[dotted stroke] (1,1.2)
    .. controls (0.92,0.88) and (0.78,0.53)
    .. (dotted join);
  \fill (dotted join) circle (0.4pt);
\end{tikzpicture}%
}}
\quad=\quad
\vcenter{\hbox{%
\begin{tikzpicture}[
  x=1cm, y=1cm, thick, baseline=0ex,
  squarednode/.style={rectangle, draw=black, minimum size=1pt}
]
  \node[squarednode] (A) at (0,0)
    {$\scriptstyle\hspace{2pt}g\hspace{2pt}$};

  \draw (A.north) -- +(0,1)
    node[xshift=5pt] {\tiny $X$};
  \draw (A.south) -- +(0,-1)
    node[xshift=3pt] {\tiny $X$};

  \draw[dotted stroke] (0,0.5)
    .. controls +(0.7,0.8) and +(1,0.2) .. (0,-0.7)
    coordinate[pos=0.78] (dotted join);

  \draw[dotted stroke] (1,1.2)
    .. controls (1.05,0.65) and (0.95,-0.20)
    .. (dotted join);
  \fill (dotted join) circle (0.4pt);
\end{tikzpicture}%
}}
\quad=\quad
\vcenter{\hbox{%
\begin{tikzpicture}[
  x=1cm, y=1cm, thick, baseline=0ex,
  squarednode/.style={rectangle, draw=black, minimum size=1pt}
]
  \node[squarednode] (A) at (0,0)
    {$\scriptstyle\hspace{2pt}g\hspace{2pt}$};

  \draw (A.north) -- +(0,1)
    node[xshift=5pt] {\tiny $X$};
  \draw (A.south) -- +(0,-1)
    node[xshift=3pt] {\tiny $X$};

  \draw[dotted stroke] (0,0.5)
    .. controls +(0.7,0.8) and +(1,0.2) .. (0,-0.7);

  \draw[dotted stroke] (1,1.2)
    .. controls (1.05,0.95) and (1.05,0.50) .. (1.05,0.15)
    .. controls (1.05,-0.20) and (0.85,-0.73) .. (0,-0.9);
  \fill (0,-0.9) circle (0.4pt);
\end{tikzpicture}%
}}\,,
\]
where the second equality of follows from the Frobenius properties, and the other are simply the consequence of naturality and associativity of $A$-action. Therefore, $E(g)$ is an $A$-module morphism, and this proves (i)

(ii) We apply (i) to $X=A$ and $g \in \Auta(A)$. Since $A$ is a simple $A$-module, $E(g)$ is a scalar multiple of $\id_A$. By taking trace in $\CC$, we find $E(g)= \tr_\CC(g) \id_A$.
On the other hand,
$$
E(g)= \begin{tikzpicture}[squarednode/.style={rectangle, draw=black,  minimum size=3pt}, baseline=0ex]
\node[squarednode](A) at (0,0) { $\hspace{2pt}\scriptstyle g \hspace{2pt}$}; 
\draw[thick, densely dotted] (A.90)-- +(0,1) node [xshift=5pt] {{\tiny $A$}};
\draw[thick, densely dotted] (A.270)-- +(0,-1)  node [xshift=3pt]{{\tiny $A$}};
\draw[thick,densely dotted] (0,.5) .. controls +(.7,.8) and +(1,.2) .. (0,-.7);
\end{tikzpicture} = 
\begin{tikzpicture}[squarednode/.style={rectangle, draw=black,  minimum size=1pt}, baseline=.5ex]
\node[squarednode](A) at (.6,.6) {$\scriptstyle g$}; 
\node[squarednode](B) at (0,.7) { $\scriptstyle g$}; 
\draw[thick, densely dotted] (B.270) -- +(0, -1.6) node [xshift=5pt] {{\tiny $A$}};
\draw[thick, densely dotted] (0,1.3) node [xshift=5pt] {{\tiny $A$}} -- (B.90) ;
\draw[thick,densely dotted] (0,-.1) .. controls +(.4,.2) and +(0, -.1) .. (A.270);
\draw[thick,densely dotted] (A.90) .. controls +(.2,1) and +(1.8,1) .. (0,-.6);
\end{tikzpicture} = 
\begin{tikzpicture}[squarednode/.style={rectangle, draw=black,  minimum size=1pt}, baseline=1ex]
\node[squarednode](A) at (.6, .6) { $\scriptstyle g$}; 
\node[squarednode](B) at (0,0.7) { $\scriptstyle g$}; 
\draw[thick, densely dotted] (0,1.4) node [xshift=5pt] {{\tiny $A$}} -- (B.90) ;
\coordinate (C) at (0.6, 0);
\draw[thick, densely dotted] (A.270) -- (C)  ;
\draw[thick, densely dotted] (B.270) -- +(0, -1.5) node [xshift=5pt] {{\tiny $A$}};
\draw[thick,densely dotted] (A.90) .. controls +(.2, .8) and +(1,.8) .. (C);
\draw[thick,densely dotted] (C) .. controls +(0,-.2) and +(-.2,-.2) .. (0,-.6);
\end{tikzpicture} = \tr_\CC(g) g\,.
$$
Therefore, $\tr_\CC(g) = 0$ if $g \ne \id_A$, and our assertion  follows.
\end{proof}

Here is our first main theorem in this paper. An explicit construction of elements 
of $\Auta(A)$ will be given in Section 7 using the canonical basis of $\CC(A,A)$.

\begin{thm} \label{t:finite}
  Let $\CC$ be a ribbon finite tensor category (not necessarily semisimple), and $A$ a commutative algebra in $\CC$ such that $A$ is simple (as an $A$-module), $\CC(\1, A)$ is 1-dimensional, and $\1$ is a direct summand of $A$. If $d(A) \ne 0$, then $\Auta(A)$ is a linear independent subset of $\CC(A,A)$, and hence $|\Auta(A)| \le \dim_\BC \CC(A,A)$. In particular, if $\CC$ is a modular fusion category and $A \in \CC$ is a condensable algebra, then $|\Auta(A)| \le \dim_\BC \CC(A,A)$.
\end{thm}
\begin{proof}
Suppose $g_1, \cdots, g_n$ are distinct elements of $\Auta(A)$ such that 
$$
0 = \gamma_1 g_1 + \cdots + \gamma_n g_n
$$
for some  scalars $\gamma_1, \cdots, \gamma_n$. Fix $i$ and  multiply $g_i^{-1}$ to the preceding equation to obtain 
$$
0 = \sum_{j=1}^n\gamma_jg_i^{-1}g_j\,.
$$
Since $g_i^{-1} g_j \ne \id_A$ for $j\ne i$, by taking trace of the equation and Lemma \ref{p5.26}, we find that $\gamma_i=0 $ as $d(A)\ne 0.$ So $g_1, \cdots, g_n$  are linearly independent.

Since $\CC(A,A)$ is a finite-dimensional $\BC$-linear space and $\Auta(A)$ is a linear independent over $\BC$,  $\Auta(A)$ must be finite and its order is bounded by $\dim_\BC \CC(A,A)$. 
\end{proof}
For a similar result in the setting of subfactors, see  \cite[Lem. 2.8]{Xu2010}. 
The group $\Auta(A)$ was previously studied in \cite{Ki2002}. In particular, it was shown that every compact subgroup of $\Auta(A)$ is finite. Thus, Theorem~\ref{t:finite} implies Kirillov's result in \cite{Ki2002} where  $A^G = \1$ is assumed.

Here is an application of Theorem \ref{t:finite} to vertex operator algebra extension. 
\begin{cor}\label{l5.4} Let $V$ be rational and $C_2$-cofinote, selfdual, simple vertex operator algebra,  and $A$ an extension vertex operator vertex operator algebra. Then $${\rm Gal}(A/V)=\{g\in{\rm Aut}(A)|g|_{V}=\id_V\}$$ is a finite group.
\end{cor}
\begin{proof} From the assumptions, the $V$-module category $\MM_V$ is a modular tesnor category \cite{huang2005vertex} and $A\in \MM_V$ is a condensable algebra \cite{HKL}. The result follows immediately from Theorem \ref{t:finite} by noting that ${\rm Gal}(A/V)={\rm Aut}_{\MM_V}^{\rm alg}(A)$.
\end{proof}
We next prove that Corollary \ref{l5.4} holds with much weeker assumptions. For this purpose
we need to recall a result from  \cite{DM1998}. Let $A$ be a vertex operator algebra and  $U,W$ be homogeneous subsapces of $A.$ For $m\in\BZ,$ set
\[
Z_m
=
\left\langle
\sum_{k=m}^{\infty}u_k v
\mid
u\in U,\ v\in W
\right\rangle 
\]
which is a subspace of $A.$
Then if $m$ is sufficiently small,
\[
\psi_m:u\otimes v\mapsto
\sum_{k=m}^{\infty}u_k v,
\qquad
u\in U,\ v\in W
\]
is a linear isomorphism (see the proof of \cite[Prop. 2.1]{DM1998}) between $U\otimes W$ and $Z_m.$ If $U,W$ are modules for some automorphism group $G$ of $A$, then $\psi_m$ is a $G$-module isomorphism.

\begin{prop} Let $A$ be a simple vertex operator algebra and $V$ simple conformal subalgebra  such that $\Hom_V(A,A)$ is finite dimensional.   Then ${\rm Gal}(A/V)$ is a finite group.
    \end{prop}

\begin{proof} Fix $g\in {\rm Gal}(A/V)$. Since
each homogeneous subspace of $A$ is finite dimensional, we have decomposition
$$A=\bigoplus_{\lambda\in \BC^*} A^{\lambda}$$
where $A^\lambda$ is the generalized eigenspace with eigenvalue $\lambda.$
From the definition of automorphism, we find
$$
(g-\ld \mu)^n Y(u,z)v= \sum_{k=0}^n \binom{n}{k} Y((g-\ld)^k \ld^{n-k} u,z)g^k(g-\mu)^{n-k} v
$$
for any nonnegative integer $n$ and $u\in A^{\lambda}, v\in A^{\mu}$ by induction on $n$. 
This implies
that $Y(u,z)v\in A^{\ld\mu}[[z,z^{-1}]].$ 
In particular, each 
$0\ne A^{\lambda}$ is an irreducible $A^1$-module \cite{DM1997}. From the assumption
$\dim\Hom_V(A,A)<\infty,$
$\{\lambda|A^{\lambda}\ne 0\}$ is a finite cyclic group generated by some $n$-th primitive root of unity. So
$A=\oplus_{r=1}^{n}A^r$ where $A^r$ is the generalized eigenspace of $g$ with eigenvalue $\lambda_r=e^{2\pi i(r-1)/n}.$ We next prove that $g$ is semisimple.

We need a result from linear algebra. Let $u$ be a generalized eigenvector of group element $h$ with eigenvalue $\lambda.$ We call the least positive integer $d$ satisfying $(h-\lambda)^du=0$ the degree of $u.$ If $u,v$ are generalized eigenvectors with degrees $d_1,d_2,$ respectively, then $u\otimes v$ is a generalized eigenvector of degree $d_1+d_2-1$ (cf. \cite{horn2012matrix}).

From the assumption we know that for each $r$ there exists a smallest positive integer $d_r$ such that $(g-\lambda_r)^{d_r}A^r=0.$ We claim that $d_r=d_1$ for all $r.$ Since $A^r$ is an irreducible $A^1$-module,
we know from \cite{DM1997} that $A^r$ is spanned by $u_nv$ for
$u\in A^1, n\in \BZ$ and a fixed  $v\in A^r\setminus\{0\}$. Let $v$ be an eigenvector. From the discussion above 
we see that $(g-\lambda_r)^{d_1}A^{r}=0$. That is, $d_r\leq d_1.$  Now we take $u\in A^1$ be a homogeneous generalized eigenvector of degree $d_1.$ Let $U$ be 
the $\langle g\rangle$-module generated by $u$ and $W=\BC v.$
Then $u\otimes v$ is a generalized eigenvector of degree $d_1$ with eigenvalue $\lambda_r.$  Thus, for small $m,$  $\psi_m(u\otimes v)\in A^r$ is a generalized eigenvector of degree $d_1$ with eigenvalue 
$\lambda_r.$ So $d_r=d_1$ for all $r.$

Now we prove $d_1=1.$  Assume that $d_1\geq 2$ and  $u\in A^1$  a generalized eigenvector of degree $d_1$.  We can similarly show that for small $m$,
$\psi_m(u\otimes u)\in A^1$ is a generalized eigenvector of degree $2d_1-1>d_1,$ a contradiction.
Thus $d_1=1$ and $g$ is semisimple. Consequently, the order of $g$ is finite.

Clearly, ${\rm Gal}(A/V)$ is a subset of $\Hom_V(A,A).$ Then ${\rm Gal}(A/V)$ acts on $ \Hom_V(A,A)$ faithfully  by left composition.  So we can regard ${\rm Gal}(A/V)$ as a subgroup of ${\rm GL}(X)$ where $X=\Hom_V(A,A)$ is finite dimensional. It follows from \cite{S1911} that every finitely generated  subgroup of ${\rm Gal}(A/V)$ is a finite group.

Finally, we show that ${\rm Gal}(A/V)$ is a finite group. From \cite{DLMcomp1996}, we know that for any finite group $G$ of ${\rm Gal}(A/V)$, 
$$A=\bigoplus_{\chi\in\irr(G)}W_{\chi}\otimes A_\chi$$
where $W_\chi$ is the irreducible $G$-module affording to the character $\chi$ and the multiplicity space $A_{\chi}$ is an irreducible $A^G$-module. Moreover, $A_{\chi}$ and $A_{\gamma}$ are inequivalent $A^G$-modules if $\chi\ne \gamma.$ This implies  that 
$$o(G)=\sum_{\chi\in\irr(G)}\chi(1)^2=\dim\Hom_{A^G}(A,A)\leq \dim\Hom_{V}(A,A)$$
as $V$ is a sub VOA of $A^G.$ So the order of any finite subgroup of  ${\rm Gal}(A/V)$ is bounded by $\dim\Hom_{V}(A,A).$ This forces  ${\rm Gal}(A/V)$ to be a finite group.
\end{proof}

We now come back to the categorical discussion.
Let $G$ be any subgroup of $\Auta(A)$. For any $x \in \CC$,  $G$ acts on $\CC(x,A)$ by composition, and so $\CC(x,A)$ can be decomposed into a direct sum isotypic components, namely
$$
\CC(x,A) =\bigoplus_{\chi\in\irr(G)} \CC(x,A)_\chi = \bigoplus_{\chi\in\irr(G)} e_\chi\CC(x,A)
$$
where $\CC(x,A)_\chi$ is the isotypic component of $\CC(x,A)$ associated with the irreducible character $\chi$ of $G$ and $e_\chi$ is the central idempotent of $\BC[G]$ corresponding to  $\chi$ which is
\[
e_\chi=\frac{\chi(1)}{|G|}
\sum_{g\in G}\chi(g^{-1})g.
\]
 Obviously $\{e_{\chi}\mid  \chi \in \irr(G)\}$ is a set of orthogonal idempotents in $\CC(A,A)$ satisfying 
$$ 
\sum_{\chi} e_{\chi} = \id_A\,.
$$
In term of the $G$-action,  $A$ can be stratified by the irreducible character $\chi$ of $G$ as 
$$
A_\chi = \sum_{x  \in \iC} \CC(x,A)_\chi\, x.
$$

Let $e_{\chi} = \iota_{\chi} \pi_{\chi}$ be the epi-monic factorization of $e_{\chi}$. Then $A_\chi$ is the image of $\pi_\chi$, $A_\chi$ is the subobject of $(A_\chi, \iota_\chi)$ of $A$, and $\pi_\chi \iota_\chi =\id_{A_{\chi}}$. For the trivial character $\chi_0$ of $G$, we simply denote by $\iota_0$ and $\pi_0$ for $\iota_{\chi_0}$ and  $\pi_{\chi_0}$ respectively, and abbreviate $A_{\chi_0}$ as $A^G$. 
\begin{rmk} \label{r:repr}
It worth noting that $\CC(x, A_\chi) \cong \CC(x, A)_\chi$ via the natural isomorphism $f \mapsto \i_\chi f$, and its inverse is given by $h \mapsto \pi_\chi h $. Thus, $A_\chi$ is the representing object of the contravariant functor $\CC(-, A)_\chi : \CC \to \Vec$. This natural isomorphism will be used throughout the remainder of this paper. 
\end{rmk}


\begin{prop} \label{p:5.16} Let $G$ be any subgroup of $\Auta(A)$. 
\begin{enumerate}
    \item For any $\chi \in \irr(G)$,
     $$
    d(A_\chi) = \frac{d(A)}{|G|} \chi(1)^2 = d(A^G) \chi(1) ^2\ne 0.
    $$
    In particular,  $|G| \mid d(A)$.
    \item The subobject $B=A^G$ of $A$ is a condensable subalgebra,
    \item $A$ is a condensable algebra in $\CBo$ with $d_B(A)=|G|$, and $G$ is subgroup of $\Aut_{\CBo}^{\alg}(A,A)$.
    \item For any $\chi \in \irr(G)$, $A_\chi$ is a local $B$-module and $d_B(A_\chi) =\chi(1)^2$. 
    \item  If $G$ is abelian, $A_\chi$ is invertible, the set $\{[A_\chi] \mid \chi \in \irr(G)\}$ of isomorphism classes  forms a group, under the product $\o_B$, and it is isomorphic the character group 
    $\irr(G)$. Moreover, $\CC_B(A,A)= \BC[G]$ and $G = \Aut_{\CBo}^{\alg}(A)$.
\end{enumerate}
\end{prop}
\begin{proof}
Since $\tr_\CC(uv) = \tr_\CC(vu)$ for any $u,v \in \CC(A,A)$, $\tr_\CC$ defines a class function on $G$. Denote the restriction of  $\tr_\CC$ to $G$ by $\tr_\CC|_G.$
By Lemma \ref{p5.26}, 
$$
\tr_\CC|_G = \frac{d(A)}{|G|} \chi^G_{reg} = \frac{d(A)}{|G|} \sum_{\chi \in \irr(G)} \chi(1)\chi
$$
where $\chi^G_{reg}$ denotes the regular character of $G$. On the other hand, for $g \in G$,
\begin{align*}
 \tr_\CC(g) &=\, \sum_{x \in \iC} \Tr_{\CC(x,A)}(g) d(x) = \sum_{x \in \iC}\sum_{\chi \in \irr(G)} N_{x,\chi}\,\chi(g) d(x) \\
 & =\, \sum_{\chi \in \irr(G)} \Big(\sum_{x \in \iC} N_{x,\chi}\, d(x)\Big) \chi(g)\,,
\end{align*}
where $N_{x, \chi}$ is the multiplicity of $\chi$ in $\CC(x,A)$ or $N_{x,\chi}\chi(1) = \dim(\CC(x,A)_\chi)$. Thus, by the linear independence of irreducible characters, we find
$$
\sum_{x \in \iC} N_{x,\chi}\, d(x) = \frac{d(A)}{|G|}\chi(1)\,.
$$
By the definition of $A_\chi$, we have
$$
d(A_\chi) = \sum_{x \in \iC} N_{x,\chi} \,d(x)\chi(1) = \frac{d(A)}{|G|}\chi(1)^2 \ne 0\,.
$$
In particular,
$$
d(A^G )= d(A_{\chi_0}) =  \frac{d(A)}{|G|}\chi_0(1)^2 = \frac{d(A)}{|G|} 
$$
 is an cyclotomic integer. Therefore, $|G|\mid d(A)$ and  this completes the proof of (i).

(ii) For any $\chi \in \irr(G)$ and $\chi' \in \hat{G}$,   $\chi\chi' \in \irr(G)$, and we let
$$
m_{\chi, \chi'}:= \pi_{\chi\chi'}m_A(\i_\chi \o \i_{\chi'}) : A_\chi \o A_{\chi'} \to A_{\chi\chi'}\,.
$$
Then 
\begin{align*}
\iota_{\chi\chi'} m_{\chi, \chi'} =& e_{\chi\chi'} m_A(\i_\chi \o \i_{\chi'})  = \frac{\chi\chi'(1)}{|G|}\sum_{g \in G} \chi\chi'(g^{-1})  m_A(g\i_\chi \o g\i_{\chi'}) \\
=& \frac{\chi(1)}{|G|}\sum_{g \in G} \chi(g^{-1}) m_A(g\i_\chi \o \i_{\chi'}) = m_A(e_\chi \i_\chi \o \i_{\chi'}) = m_A(\i_\chi \o \i_{\chi'})\,.
\end{align*}
It follows from  \cite[Lemma 5.3]{DNRX2025} that $\iota_{\chi\chi'} m_{\chi, \chi'} \ne 0$. By Remark \ref{r:repr},
\begin{equation}\label{eq:ten_prod}
  m_{\chi, \chi'}\in \CC(A_\chi \o A_{\chi'}, A_{\chi\chi'}) \setminus \{0\} \,.  
\end{equation}

We apply \eqref{eq:ten_prod} for $\chi = \chi'=\chi_0$, and let $B=A^G$.  Setting $u_B := \pi_0 u_A$ and $m_B: =\pi_0 m_A (\iota_0 \o \iota_0)$ makes $B$ to be an algebra in $\CC$ and $\iota_0: B \to A$ is an algebra monomorphism. Thus, $(B, \iota_0)$ is subalgebra of $A$, and hence a condensable algebra by Proposition \ref{p:positivity}.

(iii) Since $A$ is commutative, $A$ is a local $B$-module. Let $p_B: A \o A \to A \o_B A$ be the epimorphism defining $A \o_B A$. Then the multiplication $m_A$ of $A$ can be factorized as $m_A = \ol m_A p_B$ for a unique $B$-module morphism $\ol m_A: A \o_B A \to A$. Since
$$
\CB(B, A) \cong \CC(\1, A) \cong \BC,
$$
$A$ is a connected commutative algebra in $\CBo$ with multiplication $\ol m_A$ and unit map $\ol u_A :=\iota_0$. The morphism $\mu_A = m_A (\id_A \o \iota_0): A\otimes B\to A$ defines a $B$-module structure on $A.$
For $g \in \Auta(A)$, $g \in \CB(A,A)$ as 
$$
g \mu_A = g m_A (\id_A \o \iota_0) = m_A (g \o g \iota_0) =m_A (g \o \iota_0)=\mu_A(g\otimes \id_B)\,.
$$
One also needs to check $g \ol m_A = \ol m_A(g \o_B g)$ and $g \ol u_A=\ol u_A,$ but we leave it as an exercise to the readers. It follows directly from (i) that $d_B(A) =|G|$.\\
(iv) 
For any $\chi \in \irr(G)$,  it follows from \eqref{eq:ten_prod} that $\mu_\chi:=m_{\chi, \chi_0} \in \CC(A_\chi\o B, A_\chi) \setminus\{0\}$, and it
defines a (right) $B$-module action $\mu_\chi$ on $A_\chi$. So, $A_\chi$ is a $B$-submodule of $A$, which implies $A_\chi$ is a local $B$-module or $A_\chi \in \CBo$. By (i),
$$
d_B(A_\chi)=\frac{d(A_\chi)}{d(B)} = \chi(1)^2\,.
$$

(v) Since $G$ abelian, $\irr(G)=\hat{G}$.  Now we apply \eqref{eq:ten_prod} to the condensable algebra $A$ in $\CBo$ by (ii) and (iii), where $B=A^G$, i.e.,
$$
m_{\chi, \chi'} : A_\chi \o_B A_{\chi'} \to A_{\chi\chi'} 
$$
is a nonzero $B$-module morphism for any  $\chi, \chi'  \in \hat G$. 

 Let $M_\chi$ be any simple $B$-submodule of $A_\chi$. Then the restriction of the $B$-module morphism $m_{\chi, \ol\chi}: M_\chi \o_B M_{\ol\chi} \to B$ is still nonzero  by \cite[Lemma 5.3]{DNRX2025} again.  Thus, $M_{\ol\chi}$ and $M_\chi^*$ are isomorphic $B$-modules. Since $M_\chi$ is arbitrarily chosen, all simple $B$-submodules of $A_\chi$ are isomorphic, and so $A_\chi \cong n_\chi M_\chi$ for some $n_\chi \in \BN$. Thus
$$
1=\chi(1)^2 = d_B(A_\chi) = n_\chi d_B(M_\chi)\,.
$$
In particular, $d_B(M_\chi) \in \BZ$ which forces $d_B(M_\chi)=1=n_\chi$ for $\chi \in \irr(G)$.  Thus,  $A_\chi = M_\chi$ is a simple $B$-module. Now, the first statement follows immediately from the following Lemma \ref{l:G-grading}.  Thus, $\dim(\CBo(A,A)) = |G|$ by Schur-Weyl duality. Hence, $\CBo(A,A) = \BC[G]$ and $\Aut_{\CBo}^{\alg}(A) = G$ by Theorem \ref{t:finite}.
\end{proof}
\begin{lem}\label{l:G-grading}
Let $A$ be a  connected $G$-graded algebra in a fusion category  $\DD$ for some finite group $G$, i.e.,
$$
A = \bigoplus_{g \in G} A_{g}
$$
where $\{A_g \mid  g\in G\}$ is a set of simple objects of $\DD$ with $A_1 = \1$ such that the image of $A_g \o A_h$ under $m_A$ is $A_{gh}$ and $u_A$ is the inclusion map. Then, $\{[A_g]|g \in G\}$ is a group of distinct isomorphism classes of invertible objects under tensor product, and it is isomorphic to $G$. 
\end{lem}
\begin{proof}
Let $m_{g,h}: A_g \o A_h \to A_{gh}$ be the restriction of $m_A$ on $A_g \o A_h$. To prove that $A_g$ is invertible,  it suffices to show that $m_{g, g^{-1}}: A_g \o A_{g^{-1}}\to \1$ is a monomorphism. Since $m_{g^{-1}, g}: A_{g^{-1}} \o A_{g} \to \1$ is nonzero, $A_g^* \cong A_{g^{-1}}$. One can consider the left dual object $A_{g^{-1}}$ of $A_g$ with evaluation $m_{g^{-1}, g}$ and coevaluation $b: \1 \to A_g \o A_{g^{-1}}$.
Assume $\DD$ is strict and so the associativity of $m_A$ implies $\id_g \o m_{g^{-1}, g} = m_{g,g^{-1}} \o \id_g$, where $\id_{A_g}$ is simply abbreviated as $\id_g$.  Now we have
\begin{align*}
 b m_{g, g^{-1}}  = & b \o m_{g, g^{-1}}  = (b\o \id_g \o \id_{g^{-1}}) (\id_g \o \id_{g^{-1}} \o m_{g, g^{-1}} )\\
      =& (b\o \id_g \o \id_{g^{-1}}) (\id_g \o m_{g^{-1}, g} \o \id_{g^{-1}} ) =  \id_g \o \id_{g^{-1}}\,.
\end{align*}
Therefore, $m_{g, g^{-1}}$ is monomorphism, and so $A_g$ is invertible. In particular, $d_+(A_g) = 1$ and $m_{g,h}: A_g \o A_h \to A_{gh}$ is an isomorphism.

For $g \in G$, $d_+(A_g \o A) = |G| =d_+(A)$. Since $\DD_A(A_g \o A, A) \cong \DD(A_g, A)$, $[A_g,A]_\DD A$ is a summand of $A_g \o A$. Therefore, 
$$
[A_g,A]_\DD |G| = [A_g, A]_\DD d_+(A) \le d_+(A_g\o A) = |G|, 
$$
which implies $[A_g, A]_\DD = 1$. Thus, $A_g \not\cong A_h$ if $g \ne h$, and $\{[A_g] \mid g \in G\}$ forms a group under the tensor product. Obviously, this group is isomorphic to $G$.
\end{proof}

\section{Categorical Orbifolds} 

In this section, we develop the theory of categorical orbifolds associated with algebra automorphisms of a condensable algebra $A$. We introduce twisted $A$-modules and show that they form the homogeneous components of a $G$-graded extension of $\CAo$. Conversely, we show that any $G$-extension of $\CAo$ in $\CA$  realizes $G$ as a subgroup of $\Auta(A).$ This framework generalizes the classical orbifold construction and provides a categorical interpretation of group extensions in terms of twisted modules.

\begin{defn} Let $A$ be a condensable algebra in an MTC $\CC$, and $g \in \Auta(A)$.
    An $A$-module $M \in \CA$ is called  $g$-\emph{twisted} if the $A$-action morphism $\mu_M: M \o A \to M$ satisfies
$$
\mu_M R_{A,M} R_{M,A} = \mu_M(\id_M \o g)
$$
which is equivalent to the condition
$$
\mu_M R_{M,A}^{-1} R_{A,M}^{-1} = \mu_M(\id_M \o g^{-1})\,.
$$
We denote the full subcategory of $g$-twisted $A$-module by $\CA(g)$. By definition, $\CA(1) = \CAo$. 

A \emph{twisted $A$-module subcategory}   is a full fusion subcategory of  $\CA$  containing $\CAo$ whose objects are $g$-twisted $A$-modules for some $g \in \Auta(A)$.
\end{defn}
We first show existence of $g$-twisted module. 
\begin{lem}\label{l:nontrivial}
    For any $g \in \Auta(A)$, $\CA(g)$ is nontrivial. Moreover, for $X \in \CA(g)$
    $$
    \frac{1}{d_A(X)} \rho(X) = g
    $$
    where  $\rho: K(\CA) \to \CC(A,A)$ is the canonical surjection.
\end{lem}
\begin{proof}
    Let $H$ be the cyclic group generated by $g$, $\hat{H}$ its character group and $B=A^H$. Then $H=\Aut_{\CBo}^{\alg}(A)$ by  Proposition \ref{p:5.16} (v), $\{A_\chi\mid \chi \in \hat{H}\}$ is a group of invertible objects in $\CBo$ isomorphic to $\hat{H}$, and
    $$
A = \bigoplus_{\chi \in \hat{H}} A_\chi\,.
    $$
     We identify $A_{\chi}$ with $\chi$ for any $\chi \in \hat{H}$, and $\hat{\hat{H}}$ with $H$. The $S$-matrix of $\CBo$ determines an $\hat{\hat{H}}$-grading of $\CBo$, namely $(\CBo)_h$ is the full subcategory of $X \in \iCBo$ such that
    $$
    \chi(h) = \frac{\ol S_{X,\chi}}{d_B(X)} \,.
    $$
 Since $X \o \chi \in \iCBo$, this equality is equivalent to
\[
R_{\chi,X}R_{X,\chi}
=
\chi(h)\,\id_{X\otimes_B\chi}
\qquad
\text{for all }\chi\in H.
\]
Note that $\theta_\chi=\id_\chi$ for $\chi\in\widehat{H}$, $A=\sum_{\chi\in\widehat{H}}\chi\in (\CBo)_1$.
For any $X  \in \irr((\CBo)_g)$, $M:=X \o_B A$ is an  $A$-module in $(\CBo)_g$. Now, we show that $M$ is a $g$-twisted $A$-module by verifying the definition,
    \begin{align*}
         \mu_{M} R_{A,M}R_{M,A}  = \mu_{M} & \left(\sum_{\chi \in \hat{H}}  R_{\chi,M}R_{M,\chi} \right) \\
        &= \mu_{M} \left(\sum_{\chi \in \hat{H}} \chi(g)\id_{M \o_B \chi} \right)
     = \mu_{M}(\id_{M} \o g)\,. 
    \end{align*}
    Thus, $M \in \CA(g)$.

    For $X \in \CA(g)$,  we have
    \begin{align*}
    d(A)\rho(X) = & \begin{tikzpicture}[baseline=-1ex, scale=.55]
\begin{knot}[
end tolerance =1pt,
consider self intersections,
ignore endpoint intersections=false,
clip width=0pt,
clip radius=0pt,
]
\strand[thick] (-.5,1.3) arc[start angle=115,delta angle=-340,x radius=1, y radius=1.5];
\strand[thick, densely dotted] (-1, 2.5) 
to [out=-80, in=30] (-1.05,-.4);
\strand[thick, densely dotted] (.9, -.4) 
to [out=70, in=200] (1.2, .1)
to [out=0, in=90] (1.3, -1.8)
;
 \flipcrossings{1};
\end{knot}
\node at (-1.2,.7) { $\scriptstyle X$}; 
\end{tikzpicture} 
\!=\!
\begin{tikzpicture}[baseline=-1ex, scale=.55]
\begin{knot}[
end tolerance =1pt,
consider self intersections,
ignore endpoint intersections=false,
clip width=0pt,
clip radius=0pt,
]
\strand[thick] (-.8, 1) arc[start angle=140,delta angle=230,x radius=1, y radius=1.5];
\strand[thick] (-.6, 1.2) arc[start angle=125,delta angle=-113,x radius=1, y radius=1.5];
\strand[thick, densely dotted] (-1.2, 2.5) 
to [out=-70, in=130] (-.5,.8)
to [out=-50, in=250] (0, 1.35);
\strand[thick, densely dotted] (.05, 1.6)
to [out= 80, in=135] (.4, 1.70) 
to [out= -45, in=50] (.5, 1.38) ;
\strand[thick, densely dotted] (.43, 1.15) 
to [out=-120, in=60] (-.95,-.6) ;
\strand[thick, densely dotted] (.92, -.4) 
to [out=70, in=200] (1.2, .1)
to [out=0, in=90] (1.3, -1.8);
 \flipcrossings{1};
\end{knot}
\node at (-1.2,.7) { $\scriptstyle X$}; 
\end{tikzpicture}  
=\!
\begin{tikzpicture}[baseline=-1ex, scale=.55]
\begin{knot}[
end tolerance =1pt,
consider self intersections,
ignore endpoint intersections=false,
clip width=4pt,
clip radius=4pt,
]
\strand[thick] (0.4, 1.3) arc[start angle=70,delta angle=350,x radius=1, y radius=1.5];
\strand[thick, densely dotted] (1, 2) 
to [out=250, in=110] (.4,.8)
to [out=-70, in=110] (1.4,.6)
to [out=-70, in=70] (1.05,.2)
;
\strand[thick, densely dotted] (1, -.7) 
to [out=60, in=180] (1.3, -.2)
to [out=0, in=90] (1.3, -1.8)
;
\end{knot}
\node at (-1,.7) { $\scriptstyle X$}; 
\end{tikzpicture}
\!=\! 
\begin{tikzpicture}[squarednode/.style={rectangle, draw=black}, baseline=-1ex, scale=.55]
\node[squarednode](A) at (1.3, 1.2) { $ \scriptstyle g $}; 
\draw[thick, densely dotted] (1.3,2) -- (A.90);
\begin{knot}[
end tolerance =1pt,
consider self intersections,
ignore endpoint intersections=false,
clip width=4pt,
clip radius=4pt,
]
\strand[thick] (0,0) ellipse [x radius=1, y radius=1.5];
\strand[thick, densely dotted] (A.270) 
to [out=270, in=70] (1,.2)
;
\strand[thick, densely dotted] (.9, -.7) 
to [out=60, in=180] (1.3, -.2)
to [out=0, in=90] (1.3, -1.8)
;
\end{knot}
\node at (-1.1,.7) { $\scriptstyle X$}; 
\end{tikzpicture}\! =\!  d(X) g\,,
    \end{align*}
 where the second last equality follows from the fact that $X^* \in \CA(g^{-1})$.   This proves the last assertion.
\end{proof}

\begin{thm}\label{t2}
    Let $\CC$ be an MTC, $A \in \CC$ a condensable algebra, $G$ a subgroup of $\Auta(A)$, and $B=A^G$. Then, $(\CBo)_A(g) = \CA(g)$ for all $g \in G$, and the fusion subcategory 
    $$
    (\CBo)_A = \bigoplus_{g \in G} (\CBo)_A(g)=\bigoplus_{g \in G} \CA(g)
    $$
    of $\CA$ is a $G$-extension of $\CAo$. Moreover, $\BC [G] = \CB(A,A)$ and $G=\Aut_{\CBo}^{\alg}(A)$. 
\end{thm}
\begin{proof} (i) We first show that $\DD =\bigoplus_{g \in G} \CA(g)$ is a $G$-graded fusion subcategory of $\CA$.
To prove this claim, by Lemma \ref{l:nontrivial}, it remains to show the following two assertions for any  $g, h \in \Auta(A)$:
\begin{enumerate}    
 \item[(a)] If $g \ne h$, then $\irr(\CA(g)) \cap \irr(\CA(h)) =\emptyset$.
\item[(b)] For any $M \in \CA(g)$ and $N \in \CA(h)$, $M \o_A N \in \CA(gh)$.
\end{enumerate} 

(a) Let $M$ be a simple object in $\CA(g) \cap \CA(h)$. Then we have 
$$
\begin{tikzpicture}[squarednode/.style={rectangle, draw=black,  minimum size=1pt}, baseline=-2ex]
\node[squarednode](A) at (.7,0) { $\hspace{1pt}\scriptstyle g \hspace{1pt}$}; 
\draw[thick] (0,1.2) node [xshift=5pt] {{\tiny $M$}} -- (0, -1.2) ;
\draw[thick,densely dotted]  (A.90) -- +(0, 1)  node [xshift=5pt] {{\tiny $A$}}; 
\draw[thick,densely dotted] (A.270) .. controls +(0,-.4) and +(-.1,0) .. (0,-.8);
\end{tikzpicture}\quad  =\quad 
\begin{tikzpicture} [baseline=-2ex]
\node at (0.2,1.3) { $\scriptstyle M $}; 
\node at (1.1,1.3) { $\scriptstyle A $}; 
\begin{knot}[
clip width=6pt,
clip radius=6pt,
]
\strand[thick]  (0,1.2) to [out=down, in=up] (0, -1.2);
\strand[densely dotted, thick] (1,1.2)
to [out=270, in=90] (-.4,.3)
to [out=down, in=up] (.4,-.3)
to [out=down, in=45] (0,-.7);
\flipcrossings{2}
\end{knot}
\end{tikzpicture}\quad  =\quad 
\begin{tikzpicture}[squarednode/.style={rectangle, draw=black,  minimum size=1pt}, baseline=-2ex]
\node[squarednode](A) at (.7,0) { $\hspace{1pt}\scriptstyle h \hspace{1pt}$}; 
\draw[thick] (0,1.2) node [xshift=5pt] {{\tiny $M$}} -- (0, -1.2) ;
\draw[thick,densely dotted]  (A.90) -- +(0, 1)  node [xshift=5pt] {{\tiny $A$}}; ;
\draw[thick,densely dotted] (A.270) .. controls +(0,-.4) and +(-.1,0) .. (0,-.8);
\end{tikzpicture}
$$
which implies 
$$
\begin{tikzpicture}[squarednode/.style={rectangle, draw=black,  minimum size=1pt}, baseline=-2ex]
\node[squarednode](A) at (.7,0) { $\hspace{1pt}\scriptstyle h g^{-1} \hspace{1pt}$}; 
\draw[thick] (0,1) node [xshift=5pt] {{\tiny $M$}} -- (0, -1.2) ;
\draw[thick,densely dotted]  (A.90) -- +(0, .75)  node [xshift=5pt] {{\tiny $A$}}; ;
\draw[thick,densely dotted] (A.270) .. controls +(0,-.4) and +(-.1,0) .. (0,-.8);
\end{tikzpicture} \quad = \quad 
\begin{tikzpicture}[squarednode/.style={rectangle, draw=black,  minimum size=1pt}, baseline=-2ex]
\draw[thick] (0,1) node [xshift=5pt] {{\tiny $M$}} -- (0, -1.2) ;
\node at (.9,1) {{\tiny $A$}}; 
\draw[thick,densely dotted] (.7,1) .. controls +(0,-2) and +(-.1,0) .. (0,-.8);
\end{tikzpicture}\,.
$$
By the compositions of these morphisms to $ (\mu_M \o A) \coev_A$, we find the third equality of the following diagrammatic equations:   
$$
\tr_\CC(hg^{-1})\id_M = 
\begin{tikzpicture}[squarednode/.style={rectangle, draw=black}, baseline=1ex]
\node[squarednode](A) at (1, .8) { $\scriptstyle hg^{-1}$}; 
\coordinate (B) at (0,-.8);
\draw[thick] (0,1.4) node [xshift=5pt] {{\tiny $M$}} -- (B) node [xshift=5pt] {{\tiny $M$}};
\coordinate (C) at (.4, .2);
\draw[thick, densely dotted] (0,-.2) -- (C)  ;
\draw[thick,densely dotted] (C) .. controls +(-.1,1.9) and +(0,.2) .. (A.90);
\draw[thick,densely dotted] (C) .. controls +(.3,0) and +(0,-.2) .. (A.270);
\end{tikzpicture} \quad = \quad 
\begin{tikzpicture}[squarednode/.style={rectangle, draw=black,  minimum size=1pt}, baseline=-2ex]
\node[squarednode](A) at (.7,0) { $\scriptstyle hg^{-1}$}; 
\draw[thick] (0,1.2) node [xshift=5pt] {{\tiny $M$}} -- (0, -1.2) node [xshift=5pt] {{\tiny $M$}};
\draw[thick,densely dotted] (A.270) .. controls +(0,-.4) and +(-.1,0) .. (0,-.8);
\begin{knot}[
]
\strand[densely dotted, thick] (0,.4) 
to [out=45, in=180] (.4,.8)
to [out=0, in=90] (A.90);
\end{knot}
\end{tikzpicture} \quad = \quad 
\begin{tikzpicture}[squarednode/.style={rectangle, draw=black,  minimum size=1pt}, baseline=-2ex]
\draw[thick] (0,1.2) node [xshift=5pt] {{\tiny $M$}} -- (0, -1.2)  node [xshift=5pt] {{\tiny $M$}};
\begin{knot}[
]
\strand[densely dotted, thick] (0,.4) 
to [out=50, in=180] (.4,.8)
to [out=0, in=90] (.7,.5)
to [out=270, in=45] (0,-.8)
;
\end{knot}
\end{tikzpicture} \quad = d(A)\id_M \,.
$$
The other three equalities follow from the associativity of $\mu_M$ and that $\mu_M(\id_M\o u_A)=\id_M$. So $\tr_\CC(hg^{-1}) =d(A)$  and $h = g$ by Lemma \ref{p5.26}. 

(b) Let $p_A: M \o N \to M \o_A N$ be a cokernel of $\mu_M \o N - M \o \mu^+_N: M \o A \o N \to M \o N$ in $\CC$. Then
$p_A \mu_{M \o N} R^{-1}_{M\o N, A} R^{-1}_{A, M \o N}$ is equal to
$$
\begin{tikzpicture} [squarednode/.style={rectangle, draw=black,  minimum size=1pt}, baseline=-2ex]
\node (A) at (1.1,1.2) [above] { $\scriptstyle A $}; 
\node[squarednode](P) at (.25,-1.2) { $\hspace{3pt}\scriptstyle p_A \hspace{3pt}$}; 
\draw[thick](P.270) -- +(0, -.5) node [below] { $\scriptstyle M \o_A N$}; 
\begin{knot}[
end tolerance =1pt
clip width=9pt,
clip radius=9pt,
]
\strand[thick]  (P.135) to [out=90, in=270] +(0, 2.2) node  [above] { $\scriptstyle M $};
\strand[thick]  (P.45) to [out=90, in=270] +(0, 2.2) node  [above] { $\scriptstyle N $};;
\strand[densely dotted, thick] (A)
to [out=270, in=90] (-.4,.1)
to [out=270, in=90] (.8,-.5)
to [out=down, in=45] (.45,-.7);
\flipcrossings{1}
\flipcrossings{3}
\end{knot}
\end{tikzpicture} \, = \,
\begin{tikzpicture} [squarednode/.style={rectangle, draw=black,  minimum size=1pt}, baseline=-2ex]
\node (A) at (1.1,1.2) [above] { $\scriptstyle A $}; 
\node[squarednode](P) at (.25,-1.2) { $\hspace{3pt}\scriptstyle p_A \hspace{3pt}$}; 
\draw[thick] (P.270) -- +(0, -.5) node [below] { $\scriptstyle M \o_A N$}; 
\begin{knot}[
end tolerance =1pt
clip width=9pt,
clip radius=9pt,
]
\strand[thick]  (P.135) to [out=90, in=270] +(0, 2.2) node  [above] { $\scriptstyle M $};
\strand[thick]  (P.45) to [out=90, in=270] +(0, 2.2) node  [above] { $\scriptstyle N $};;
\strand[densely dotted, thick] (A)
to [out=270, in=90] (-.4,.1)
to [out=280, in=90] (.9,-.4)
to [out=down, in=45] (0.05,-.7);
\flipcrossings{3}
\flipcrossings{1}
\end{knot}
\end{tikzpicture}\, = \,
\begin{tikzpicture} [squarednode/.style={rectangle, draw=black,  minimum size=1pt}, baseline=-2ex]
\node (A) at (1.1,1.2) [above] { $\scriptstyle A $}; 
\node[squarednode](P) at (.25,-1.2) { $\hspace{3pt}\scriptstyle p_A \hspace{3pt}$}; 
\node[squarednode](g) at (1.1, .4) { $\hspace{1pt}\scriptstyle g^{-1} \hspace{1pt}$};
\draw[thick, densely dotted] (A) -- (g.90); 
\draw[thick] (P.270) -- +(0, -.5) node [below] { $\scriptstyle M \o_A N$}; 
\begin{knot}[
line join=round,
end tolerance =1pt
clip width=5pt,
clip radius=5pt,
]
\strand[thick]  (P.135) to [out=90, in=270] +(0, 2.2) node  [above] { $\scriptstyle M $};
\strand[thick]  (P.45) to [out=90, in=270] +(0, 2.2) node  [above] { $\scriptstyle N $};;
\strand[densely dotted, thick] (g.270)
to [out=270, in=65] (0.05,-.4);
\flipcrossings{1}
\end{knot}
\end{tikzpicture}
\, = \,
\begin{tikzpicture} [squarednode/.style={rectangle, draw=black,  minimum size=1pt}, baseline=-2ex]
\node (A) at (1.1,1.2) [above] { $\scriptstyle A $}; 
\node[squarednode](P) at (.25,-1.2) { $\hspace{3pt}\scriptstyle p_A \hspace{3pt}$}; 
\node[squarednode](g) at (1.1, .6) { $\hspace{1pt}\scriptstyle g^{-1} \hspace{1pt}$};
\draw[thick, densely dotted] (A) -- (g.90); 
\draw[thick] (P.270) -- +(0, -.5) node [below] { $\scriptstyle M \o_A N$}; 
\begin{knot}[
line join=round,
end tolerance=0pt,
clip width=5pt,
clip radius=5pt,
]
\strand[thick]  (P.135) to [out=90, in=270] +(0, 2.2) node  [above] { $\scriptstyle M $};
\strand[thick]  (P.45) to [out=90, in=270] +(0, 2.2) node  [above] { $\scriptstyle N $};;
\strand[densely dotted, thick] (g.270)
to [out=270, in=90] (0.2, -.1)
to [out=270, in=90] (0.75, -.6)
to [out=270, in=45] (0.5, -.75);
\flipcrossings{1}
\end{knot}
\end{tikzpicture} \, = \,
\begin{tikzpicture} [squarednode/.style={rectangle, draw=black}, baseline=-2ex]
\node (A) at (1.1,1.2) [above] { $\scriptstyle A $}; 
\node[squarednode](P) at (.25,-1.2) { $\hspace{3pt}\scriptstyle p_A \hspace{3pt}$}; 
\node[squarednode](g) at (1.1, .8) { $\hspace{1pt}\scriptstyle g^{-1} \hspace{1pt}$};
\node[squarednode](h) at (1.1, .1) { $\hspace{1pt}\scriptstyle h^{-1} \hspace{1pt}$};
\draw[thick, densely dotted] (A) -- (g.90); 
\draw[thick,densely dotted] (g.270) -- (h.90); 
\draw[thick] (P.270) -- +(0, -.5) node [below] { $\scriptstyle M \o_A N$}; 
\begin{knot}[
line join=round,
ignore endpoint intersections=true,
end tolerance=0pt,
clip width=5pt,
clip radius=5pt,
]
\strand[thick]  (P.135) to [out=90, in=270] +(0, 2.2) node  [above] { $\scriptstyle M $};
\strand[thick]  (P.45) to [out=90, in=270] +(0, 2.2) node  [above] { $\scriptstyle N $};;
\strand[densely dotted, thick] (h.270) to [out=270, in=45] (0.47, -.7);
\flipcrossings{2}
\end{knot}
\end{tikzpicture}
$$
which means 
$$
\mu_{M \o_A N} R_{M \o_A N, A}^{-1} R_{A, M \o_A N}^{-1} = \mu_{M \o_A N}(\id_M \o h^{-1} g^{-1}) = \mu_{M \o_A N}(\id_M \o (gh)^{-1})\,.
$$
So, $M \o_A N \in \CA(gh)$. This completes the proof of our claim (i).

\noindent (ii) By Proposition \ref{p:5.16} (iii), $A$ is a condensable algebra in $\CBo$ and $G$ is a subgroup of $\Aut_{\CBo}^{\alg}(A)$. For any $g \in G$, it is obvious that  $(\CBo)_A(g)$ is a full subcategory of $\CA(g)$. It follows from (i) that 
$$
\dim ((\CBo)_A(g)) = \dim((\CBo)_A^0)  =\dim (\CAo)= \dim (\CA(g))
$$
and so $(\CBo)_A(g) = \CA(g)$. Thus, $\DD$ is a fusion subcategory of $(\CBo)_A$.

 The global dimension of $\DD$ is $\dim(\CAo)|G| =\dim(\CAo)d_B(A)$ by Proposition \ref{p:5.16} (iii), and
$$
\dim((\CBo)_A) = \frac{\dim(\CBo)}{d_B(A)}=\frac{\dim(\CC)}{d(B) d(A)} = \dim(\CAo) d_B(A) =\dim(\DD)\,.
$$
So $\DD = (\CBo)_A$. \\ \\
(iii) By Schur-Weyl duality, $e_\1 K(\DD)\cong \CB(A,A)$ where 
$$
e_\1=\frac{1}{\dim(\CAo)}\sum_{X \in \iCAo} d_A(X) X\,.
$$
So $\dim_{\BC} \CB(A,A) = \dim_\BC e_\1 K(\DD) = |G|$ by Lemma \ref{l:nontrivial}. Now the last assertion follows from Theorem \ref{t:finite}.
\end{proof}

It is worth to point out  that (i) (b) of Theorem \ref{t2} has been proved previously in the VOA setting  \cite{McRae2021}, and its argument can be adapted to the categorical setting. 

In \cite[Thm. 6.1]{ENO2010}, it has been proved  that the homogeneous components of any every $G$-extension of a fusion category $\DD_1$ are invertible $\DD_1$-bimodules and hence by \cite[Thm. 7.7]{ENO2010} there is a morphism of 2-groups from $G$ to $\underline{\underline{\rm BrPic}}(\DD_1)$. By  Theorem \ref{t2}, we can conclude  the following corollary.
\begin{cor} \label{c:Pic1}
    Let $A$ be a condensable algebra in a modular fusion category $\CC$. Then there exists a morphism of 2-groups $c: \Auta(A) \to \underline{\underline{\rm BrPic}}(\CAo)$. \qedhere
\end{cor}

The following result strengthens Proposition \ref{p:5.16} (iv). 

\begin{cor} \label{c:degree} Using the same notations as in Proposition \ref{p:5.16}, for any $\chi\in \irr(G)$ and simple $B$-submodule $M_\chi$ of $A_\chi$,
$$
A_{\chi}=\chi(1)M_\chi 
$$
and $d_B(M_\chi)=\chi(1)$.
\end{cor}
\begin{proof} From  \cite[Theorem I]{DNRX2025} (or Theorem \ref{t:SW-fusion} in this paper) and Theorem \ref{t2}, for any $\chi \in\irr(G)$,
$A_{\chi}=\CBo(M_\chi,A)_\chi\otimes M_{\chi}$ for any simple $B$-submodule $M_\chi$ of $A_\chi$, and $\CBo(M_\chi,A)_\chi$ is an irreducible representation of $G$ affording the character $\chi$. Thus, $d_B(M_\chi) =\chi(1)$ by Proposition \ref{p:5.16} (iv).
\end{proof}

The first part of the following theorem was given in \cite[Prop. 5.8]{DNRX2025} already with assumption that $\CC$ is pseudounitary.
 \begin{thm}\label{p5.17}
Let $A$ be a condensable algebra in a modular fusion category $\CC$. Suppose $\BB$ is a fusion subcategory of $\CA$ and is a $G$-extension of $\CAo$ for some finite group $G$. Then the assignment $$g \mapsto \ol g := \frac{1}{\dim(\CAo)} \sum_{X \in \irr(\BB_g)}d_A(X) X$$ defines a group monomorphism from $G$ to $\KA$, and $\rho(\ol g) \in \Auta(A)$. In particular, $\rho(\ol G)$ is a subgroup of $\Auta(A)$  isomorphic to $G$. Moreover,
$$
\BB_g = \CA(\rho(\ol g))\,.
$$
\end{thm}

\begin{proof}
It follows from the proof of \cite[Thm. 3.5.2]{EGNO} that $\ol{g}\ol{h} = \ol{gh}$ for any $g, h \in G$. Obviously, $\ol g = \ol 1$ if and only if $g=1$. Thus, $\ol G$ forms a group under the multiplication of $\KA$ is isomorphic to $G$. Since $X \cdot u_A = d_A(X) u_A$ for $X \in \iCA$, $\ol g \cdot u_A = u_A$ or $\rho(\ol g)\circ u_A = u_A$. To show that $\rho(\ol g) \in \Auta(A)$, it suffices to prove that $\rho(\ol{g}) m_A = m_A (\rho(\ol{g}) \o  \rho(\ol{g}))$. 

For $x \in \iC$, let  $\{f_{x,i}\}_i$ be any basis of $\CC(x, A)$, and  $\{f^i_{x}\}_i$ be its dual basis for $\CC(A,x)$ with respect to composition, i.e.  $f_x^j \circ f_{x, i} = \delta_{i,j} \id_x$. Then 
$\id_A = \sum_{x \in \iC} \sum_{i=1}^{N_x} f_{x, i} \circ  f_x^i $, where $N_x =[x, A]_\CC$, and so 
$$
\id_A \o \id_A = \sum_{x, y \in \iC} \sum_{i=1}^{N_x} \sum_{j=1}^{N_y}   (f_{x, i} \o f_{y,j}) \circ (f_x^i \o f_y^j) \,.
$$
Thus, 
$$
\rho(\ol g) \o \rho(\ol g) = \sum_{x, y \in \iC} \sum_{i=1}^{N_x} \sum_{j=1}^{N_y} (\ol g f_{x,i} \o \,\ol g f_{y,j}) \circ (f_x^i \o f_y^j).
$$
Here, $\ol g f$ simply abbreviate  $\ol g \cdot f$ or $\rho(\ol g)f$. Since $\ol g$ has finite order, it would be convenient to assume $\{f_{x,i}\}_i$ is an eigenbasis for the $\ol g$ action on $\CC(x,A)$. Therefore, it suffices to show that 
\begin{equation}\label{eq:g_inv}
    \ol g m_A(f_{x} \o f_{y}) = \ld_{x} \ld_{y}  m_A(f_{x} \o f_{y}) = m_A( \ol g f_{x} \o\, \ol g f_{y})
\end{equation}
for  $f_x \in \CC(x,A)$ and $f_y \in \CC(y,A)$ such that $\ol g f_x = \ld_x f_x$ and $\ol g f_y = \ld_y f_y$ for some scalars $\ld_x, \ld_y$. 
However,  \eqref{eq:g_inv} can be proved with the same idea as \cite[Lemma 5.4]{DNRX2025} by using the $\KA$-module isomorphism $\psi: \CC(x, A) \to \CA(\a(x), \1)$. Then, $\ol g \psi(f_{x}) = \ld_{x} \psi(f_{x})$, and so  $\langle h_x, \ol g \psi(f_{x})  \rangle = \ld_{x} \langle h_x, \psi(f_{x})\rangle$ for any $h_x \in \CC_A(\1,\alpha(x))$. Let $X \in \irr(\BB_g)$. Then $\ol 1 X =d_A(X) \ol g$, and so the action of $\ol g$ is the same as $\frac{1}{d_A(X)} X$. In terms of diagram, we have
    \begin{equation*}
    \langle h_x, \ol g \psi(f_{x})\rangle =\frac{1}{d_A(X)}
{\color{blue}
\vcenter{\xy  0;/r1pc/:
(0,0)="o";(-1.4,1)**\dir{}, "o",
{\ellipse(1.5,1.8)^,=:a(345){-}},
(0,0)*=<20pt,10pt>[F]{\scriptstyle \psi(f_{x})},
{\ar@{-}@`{(-.1,1.5), (-2,1.1)} (0,.4);(-1.5,2.5)*=<10pt,10pt>[F]{\scriptstyle h_x}},
(-1.8,-1)*{\color{black}\scriptstyle X},
\endxy}} \quad = \ld_{x}\!\!
{\color{blue}
\vcenter{\xy  0;/r1pc/:
(0,0)="o";(-1.4,1)**\dir{}, "o",
(0,0)*=<20pt,10pt>[F]{\scriptstyle \psi(f_x)},
{\ar@{-} (0,.4);(0,2.5)*=<10pt,10pt>[F]{\scriptstyle h_x}},
\endxy}}\quad= \ld_{x}  \langle h_x, \psi(f_{x})\rangle 
\end{equation*}
and so
\begin{equation}\label{eq:half-braid}
{\color{blue}
\vcenter{\xy  0;/r1pc/:
(0,0)="o";(-1.4,1)**\dir{}, "o",
(0,-.5)*=<20pt,10pt>[F]{\scriptstyle \psi(f_{x})},
{\ar@{-}@`{(-.1,1.5), (-2,1.1)} (0,-.1);(-1.5,2.5)*=<10pt,10pt>[F]{\scriptstyle h_x}},
{\ar@{-}@`{(1,1), (-2,1.5)}|(.55){\hole} (1,2.5)*+{\color{black}\scriptstyle X};(-2,-.9)},
\endxy}} \quad = 
\ld_{x} \langle h_x, \psi(f_{x})\rangle 
{\color{blue}
\vcenter{\xy  0;/r1pc/: 
{\ar@{-} (0,2)*+ {\color{black}\scriptstyle X};(0,-1.5)},
\endxy}} \,.
\end{equation}
Here, the blue diagrams indicate they are considered in the fusion category $\CA$, and the crossing denotes the half-braiding of $\a^+(x)$. \eqref{eq:half-braid} implies
\begin{align*}
& \langle h_x \o_A h_y, \ol g(\psi(f_x) \o_A \psi(f_y))\rangle =\frac{1}{d_A(X)}
{\color{blue}
\vcenter{\xy  0;/r1pc/:
(0,0)="o";(-1.4,1)**\dir{}, "o",
{\ellipse(2.5,2.5)^,=:a(320){-}},
(-1,-.5)*=<20pt,10pt>[F]{\scriptstyle \psi(f_{x})},
(1,-.5)*=<20pt,10pt>[F]{\scriptstyle \psi(f_{y})},
{\ar@{-}@`{(-1.1,1.5), (-2.5,1.1)} (-1,-.1);(-2.5,3.5)*=<10pt,10pt>[F]{\scriptstyle h_x}},
{\ar@{-}@`{(1,1), (-1.2,1)} (1,-.1);(-1.2,3.5)*=<10pt,10pt>[F]{\scriptstyle h_y}},
(-2.2,-2)*{\color{black}\scriptstyle X},
(.1,.1)*{\ellipse(2.5,2.5)^,=:a(19){-}},
\endxy}} \\
=&\, \frac{\ld_x \langle h_x, \psi(f_{x})\rangle }{d_A(X)}
{\color{blue}
\vcenter{\xy  0;/r1pc/:
(0,0)="o";(-1.4,1)**\dir{}, "o",
{\ellipse(1.5,1.8)^,=:a(345){-}},
(0,0)*=<20pt,10pt>[F]{\scriptstyle \psi(f_2)},
{\ar@{-}@`{(-.1,1.5), (-2,1.1)} (0,.4);(-1.5,2.5)*=<10pt,10pt>[F]{\scriptstyle h_y}},
(-1.8,-1)*{\color{black}\scriptstyle X},
\endxy}} \quad
=\, \ld_x \ld_y \langle h_x, \psi(f_x)\rangle \langle h_y, \psi(f_y)\rangle\,.  
\end{align*}
Thus,
$$\ol g (\psi(f_x) \o_A \psi(f_y)) = \ld_x \ld_y \psi(f_x) \o_A \psi(f_y)\,.$$
By \cite[Lemma 5.3]{DNRX2025}, 
$$
\ol g m_A(f_x \o f_y) = \ld_x \ld_y m_A(f_x \o f_y) = m_A(\ol g f_x \o \ol g f_y)
$$
and this completes the proof of \eqref{eq:g_inv} and hence the first statement.

Since $\rho(\ol g) \in \Auta(A)$,  we can consider the $\rho(\ol g)$-twisted $A$-modules, and  let 
$$
\tilde g = \frac{1}{\dim(\CAo)}\sum_{X \in \irr(\CA(\rho(\ol g)))} d_A(X) X \in \KA\,.
$$
By Lemma \ref{l:nontrivial}, $\rho(\tilde g) = \rho(\ol g)$. Since $X\tilde g =d_A(X) \tilde g$ and $X\ol g =d_A(X) \ol g$ for $X \in \CAo$, by Theorem Theorem \ref{t:SW-fusion} (i), we find
$$
\tilde g = e_\1 \tilde g = e_\1 \ol g = \ol g \,.
$$
Thus $\irr(B_g) = \irr(\CA(\rho(\ol g)))$ which implies the last assertion.
\end{proof}
The following corollary follows immediately from Theorems \ref{t2} and  \ref{p5.17}.
\begin{cor}
    Let $\BB$ be a fusion subcategory of $\CA$. Then $\BB$ is a $G$-graded extension of $\CAo$ if, and only if, there exists a  group monomorphism $\tau: G \to \Auta(A)$ such that $\BB = \sum_{g \in G} \CA(\tau(g))$. \hfill$\square$
\end{cor}

\section{Group extensions of fusion categories}

In this section, we study group extensions of fusion categories through the bimodule decomposition associated with a fusion subcategory. Let $\DD$ be any fusion category and $\DD_1$ a fusion subcategory of $\DD$. We investigate the indecomposable $\DD_1$-bimodule decomposition of $\DD$, and obtain a canonical basis for the algebra $g_1 K(\DD) g_1$, where $g_1$ is the idempotent of $K(\DD)$ determined by $\DD_1$. Moreover, this basis defines a table algebra $\HH(\DD_1, \DD)$ and we characterize the invertible basis elements that determine the corresponding group structure. Applications of these results to 
$\CAo$ and $\CA$ will be discussed in Section 7.  The readers are referred to \cite{AB, AFM} for the definition of a table algebra which is not necessarily commutative.

These canonical basis elements are defined similarly as those in Theorem \ref{p5.17} but with respect to the Frobenius-Perron dimensions $d_+$ instead. Under the condition that $\DD_1$ \emph{pseudocentral} subcategory of $\DD$, those basis elements with Frobenius-dimension 1 form group under the multiplication $K(\DD)$.

\begin{prop}\label{p:group}
    Let $\DD$ be a fusion category, $\DD_1$  a fusion subcategory of $\DD$, and $\DD = \bigoplus_{i=1}^r \DD_i$ the indecomposable $\DD_1$-bimodules decomposition of $\DD$. Define 
\begin{equation} \label{eq:g_i}
   g_i:=\frac{1}{d_+(\DD_1)} \sum_{M \in \irr(\DD_i)} d_+(M) M \in K(\DD_i) 
\end{equation}
where $d_+$ denotes the Frobenius-Perron dimension and $d_+(\DD_i)= \sum_{M \in \irr(\DD_i)} d_+(M)^2$.
Then:
\begin{enumerate}
    \item $g_1 K(\DD) g_1$  is 1-dimensional and $g_i = \frac{d_+(\DD_i)}{d_+(M)d_+(\DD_1)}g_1 M g_1 $ for any $M \in \irr(\DD_i)$.
    \item The set $B= \{g_i \mid i=1, \dots, r\}$ is a basis of  $g_1 \KD g_1$.
    \item For any $i$, the full subcategory $\DD_{i^*}$ of $\DD$ consisting of the objects $X^*$ for all $X\in \DD_i$ is an indecomposable $\DD_1$-submodule of $\DD$. Moreover,
    \begin{equation} \label{eq:prod0}
    g_i  g_j  = \delta_{i,j^*} \frac{d_+(\DD_j)}{d_+(\DD_1)} g_1 + \sum_{k =2}^n C_{i,j}^k g_k
    \end{equation}
    where $C_{i,j}^k \in \BRn$  for $k>1$, and $d_+(g_i) = C_{i,i^*}^1 =\frac{d_+(\DD_j)}{d_+(\DD_1)} \ge 1$ for all $i$. 
    \item If, in addition, $g_1$ is in the center of $K(\DD)$, then  for any $i,j=1, \dots, r$, we have
    \begin{equation} \label{eq:prod1}
        g_i g_j =  \sum_{k=1}^r \frac{d_+(\DD_i)d_+(\DD_j)}{d_+(\DD_k)d_+(\DD_1)}\left(\sum_{Z \in \irr(\DD_k)}N_{M_i, M_j}^Z \frac{d_+(Z)}{d_+(M_i) d_+(M_j)}\right) g_k
    \end{equation}
    where $M_i \in \irr(\DD_i)$ can be arbitrarily chosen for each $i \ge 1$.  Moreover, 
    for any $M \in \irr(\DD_j)$, 
    \begin{equation} \label{eq:prod2}
    \frac{d_+(M)^2 d_+(\DD_1)}{d_+(\DD_j)} = \sum_{Z \in \irr(\DD_1)} N_{M,M^*}^Z d_+(Z)\,.
  \end{equation}
Furthermore, for any $i,j,k$, we have
\begin{equation} \label{eq:ineq1}
    C_{i,j}^k \le \frac{d_+(\DD_i)d_+(\DD_j)}{d_+(\DD_k) d_+(\DD_1)}\,. 
\end{equation}
\end{enumerate}
\end{prop}
\begin{proof}
(i) Note that  
$$
R_\DD = \sum_{X \in \iD} d_+(X)X = d_+(\DD_1)\sum_i g_i, 
$$
is a scalar multiple of the central idempotent of $K(\DD)$ corresponding to the character $d_+$ \cite[Prop. 3.3.11]{EGNO}. Thus, $
Y R_\DD =  d_+(Y) R_\DD = R_\DD Y$
for any $Y \in \irr(\DD_1)$, which implies  ${d_+}(Y) g_i = Y g_i = g_i Y$. Since $d_+(g_1)=1$, $g_1  g_i = g_i = g_i g_1$ and so $g_i \in g_1 K(\DD) g_1$ for all $i$.

 For any $i$ and any $M \in \irr(\DD_i)$, $g_1 M g_1 \ne 0$. Moreover, 
 $$
 Yg_1 M g_1  = d_+(Y) g_1 M g_1 = g_1 M g_1 Y\quad \text{for any }Y  \in \DD.
$$
  It follows from \cite[Props. 7.72 and 3.44]{EGNO} that $g_1 M g_1 = c g_i$ for some $c \in \BR^\times$. Thus $g_1 K(\DD_i) g_1$ is spanned by $g_i$.

Now, for any $M, N \in \irr(\DD_i)$, $g_1 M g_1 = c g_1 N g_1$ for some scalar $c$. Again ${d_+}(g_1)=1$ implies
$$
{d_+}(M) = c {d_+}(N) \quad \text{ or } \quad c = \frac{{d_+}(M)}{{d_+}(N)}\,.
$$
Thus, $\frac{1}{{d_+}(M)} g_1 M g_1 = \frac{1}{{d_+}(N)} g_1 N g_1$ for any $M, N \in \irr(\DD_i)$. So for any $M \in \irr(\DD_i)$,

\begin{align*}
   g_i = &\, g_1 g_i g_1  =  \frac{1}{{d_+}(\DD_1)}\sum_{N \in \irr(\DD_i)} {d_+}(N^*) g_1 N g_1 \\ 
   = &\,
   \frac{1}{{d_+}(\DD_1)}\sum_{N \in \irr(\DD_i)} {d_+}(N^*) \frac{{d_+}(N)}{{d_+}(M)} g_1 M g_1
=  \frac{{d_+}(\DD_i)}{{d_+}(\DD_1) {d_+}(M)} g_1 M g_1\,. 
\end{align*}
(ii) Obviously, $B$ is linearly independent and  $B$ spans $g_1K(\DD)g_1$ by (i).\\
(iii) Since $\DD_i$ is an indecomposable $\DD_1$-bimodule, so is $\DD_{i^*}$. Obviously, $\DD_1$ is closed under duality as $1^*=1$.

 Recall that the multiplicity $[X,Y]_\DD := \dim_\BC\DD(X,Y)$ can be extended to a bilinear form $[\cdot, \cdot]$ on $K(\DD)$ with $\irr(\DD)$ an orthonormal basis. Similarly,  duality $*$ of $\DD$ can be extended to a $\BC$-algebra antiautomorphism on $K(\DD)$, which satisfies
$$
[xz,y] = [x, yz^*] = [z, x^* y] \quad \text{for any }x,y, x \in K(\DD)\,.
$$
Under this convention, $g_i^* = g_{i^*}$ for any $i$. Thus, we have $[g_i, g_j] = \delta_{i,j}\frac{{d_+}(\DD_i)}{{d_+}(\DD_1)^2}$. By (ii), for any $i,j$, we find
\begin{equation} \label{eq:cijk}
g_i g_j = \sum_{k=1}^r C_{i,j}^k g_k
\end{equation}
for some scalars $C_{i,j}^k$. Thus, 
\begin{equation}\label{eq:cijk1}
C_{i,j}^k = [g_i g_j, g_k] \frac{{d_+}(\DD_1)^2} {{d_+}(\DD_k)}\,. 
\end{equation}
Obviously, $[g_i g_j, g_k]$ is a nonnegative real number, and so is $C_{i,j}^k$. Moreover,
\begin{equation} \label{eq:cijk2}
C_{i,j}^1 = [g_i g_j, g_1] \frac{{d_+}(\DD_1)^2} {{d_+}(\DD_1)} =  [g_i, g_{j^*}] {d_+}(\DD_1) = \delta_{i,j^*}\frac{{d_+}(\DD_i)}{{d_+}(\DD_1)}\,.
\end{equation}
Thus, 
$$
d_+(g_i) = \frac{1}{d_+(\DD_1)}\sum_{X \in \irr(\DD_i)} d_+(X^*)d_+(X) =\frac{d_+(\DD_i)}{d_+(\DD_1)} = C_{i,i^*}^1\,
$$
and  \eqref{eq:cijk} becomes
$$
g_i g_j = \delta_{i, j^*} \frac{{d_+}(\DD_i)}{{d_+}(\DD_1)} g_1+ \sum_{k=2}^r C_{i,j}^k g_k\,.
$$
Applying $d_+$ to the equation with $j=i^*$, we find
$$
d_+(g_i)^2 = d_+(g_i g_i^*) =  \frac{{d_+}(\DD_i)}{{d_+}(\DD_1)} + \sum_{k=2}^r C_{i,j}^k d_+(g_k) \ge \frac{{d_+}(\DD_i)}{{d_+}(\DD_1)}=d_+(g_i) \,.
$$
Since $d_+(g_i) >0$, $d_+(g_i) \ge 1$.\\
(iv) Now, we assume $g_1$ is in the center of $K(\DD)$. For any integer $k \ge 1$, we arbitrarily  choose an object $M_k \in \irr(\DD_k)$. By (i), we have
    \begin{align*}
        [g_{i} g_j, g_k] = &\frac{{d_+}(\DD_i){d_+}(\DD_j)[g_1 M_i M_j, g_k]}{{d_+}(\DD_1)^2{d_+}(M_i){d_+}(M_j)} 
        =  \frac{{d_+}(\DD_i){d_+}(\DD_j)[ M_i M_j, g_k]}{{d_+}(\DD_1)^2{d_+}(M_i){d_+}(M_j)} \\
        = &\, \frac{{d_+}(\DD_i){d_+}(\DD_j)}{{d_+}(\DD_1)^3{d_+}(M_i){d_+}(M_j)} \sum_{Z \in \irr(\DD_k)}N_{M_i, M_j}^Z {d_+}(Z)\,.
    \end{align*}
 Thus, by \eqref{eq:cijk1},
 \begin{equation} \label{eq:cijk3}
 C_{i,j}^k =\frac{{d_+}(\DD_i){d_+}(\DD_j)}{{d_+}(\DD_k){d_+}(\DD_1)}\sum_{Z \in \irr(\DD_k)} N_{M_i, M_j}^Z \frac{{d_+}(Z)}{{d_+}(M_i){d_+}(M_j)}\,.
\end{equation}
For any $M \in \irr(\DD_i)$, $M^* \in \irr(\DD_{i^*})$.  By \eqref{eq:cijk2}, we find 
 $$\frac{{d_+}(\DD_i)}{{d_+}(\DD_1)} =C_{i, i^*}^1 = \frac{{d_+}(\DD_i)^2}{{d_+}(\DD_1)^2|{d_+}(M)|^2}\sum_{Z\in\irr(\DD_1)}N_{M,M^*}^Z {d_+}(Z)
 $$
 which immediately implies the identity
 $$
 \frac{{d_+}(\DD_1){d_+}(M)^2}{{d_+}(\DD_i)} = \sum_{Z\in\irr(\DD_1)}N_{M,M^*}^Z {d_+}(Z)\,. 
 $$
 The  sum on right hand side of \eqref{eq:cijk3} is 
$$
\sum_{Z \in \irr(\DD_k)} N_{M_i, M_j}^Z \frac{d_+(Z)}{d_+(M_i)d_+(M_j)} \le \sum_{Z \in \irr(\DD)} N_{M_i, M_j}^Z \frac{d_+(Z)}{d_+(M_i)d_+(M_j)} = 1,
$$
which implies  inequality  \eqref{eq:ineq1}.
\end{proof}

\begin{rmk}\label{r6.2} Proposition \ref{p:group} (i) is generally false if the Frobenius-Perron dimension $d_+$ is replaced by other linear characters of $K(\DD)$.
\end{rmk}
The algebra $g_1 \KD g_1$ with the distinguished basis $B$ studied Proposition \ref{p:group} is a \emph{table algebra}, which is defined as follows  (cf. \cite{Bl}).
\begin{defn} \label{d:RBA_TA}
\begin{enumerate}
\item A \emph{reality-based algebra}, abbreviated as RBA, is a pair $(\HH, B)$ in which $\HH$ is a finite-dimensional $\BC$-algebra with a distinguished basis $B$ containing $1_\HH$ and an involutive algebra antiautomorphism $*$ such that
\begin{enumerate}
    \item $B^* = B$;
    \item the structure coefficients $C_{a,b}^c$ associated with $B$, defined by
    $$
    ab = \sum_{c \in B} C_{a,b}^c c,
    $$
    are in $\BR$ for any $a,b,c \in B$
    \item  $C_{a,b}^1 > 0$ if and only if $b=a^*$.
\end{enumerate}
    \item A \emph{table algebra} is an RBA $(\HH,B)$ such that the structure coefficients $C_{a,b}^c \in \BRn$ for any $a,b,c \in B$, where $\BR_{\ge 0}$ denote the set of all nonnegative real numbers.
\item A table algebra $(\HH, B)$ is called \emph{standard} if the Frobenius-Perron eigenvalue $\delta(a)$ of the left multiplication matrix by $a$ coincide with $C_{a,a^*}^1$ for all $a \in B$.
\item An element $b \in B$ of table algebra $(\HH, B)$ is called \emph{linear} (or \emph{thin}) if $bb^*\in \BR_{>0}1_\HH$, and the set of linear elements is denoted by $L(\HH,B)$ (cf. \cite[Def. 3.14]{Bl}).
\item Let $(\HH_1, B_1)$, $(\HH_2, B_2)$ be table algebras. An algebra homomorphism $\phi: \HH_1 \to \HH_2$ is called a homomorphism of table algebras if $\phi(B_1) \subset \BR_{> 0} B_2$, the set of all positive multiples of the elements of $B_2$. It is called epimorphism, monomorphism or isomorphism as it is such as algebra homomorphism (cf. \cite[Def. 3.11]{Bl}).
\end{enumerate}
\end{defn}

\begin{propdef}\label{pd:table-alg} Let $\DD$ be a fusion category, $\DD_1$ a fusion subcategory of $\DD$ and $\I$ the set of indecomposable $\DD_1$-bimodules in $\DD$.  For any $\MM \in \I$, we  write
   $$
    g_\MM := \frac{1}{d_+(\MM)} \sum_{M\in \irr(\MM)} d_+(M^*)M \in K(\DD)\,,
    $$
    where $d_+(\MM):=d_+(g_\MM)$. Let $B= \{g_\MM \mid \MM \in \I\}$ and  $g_1:=g_{\DD_1}$. Then  $\HH(\DD_1, \DD):=(g_1 \KD g_1, B)$ is a standard table algebra. For any $g_\MM \in B$, the following statement are equivalent:
    \begin{enumerate}
      \item $g_\MM$ is linear;
      \item $d_+(g_\MM) = 1$;
      \item $g_\MM g_{\MM^*} = g_1$.
    \end{enumerate}
     An element $g_\MM \in B$ is said to be \emph{invertible} if it satisfies any of the above conditions. The set of invertible elements of $\HH(\DD_1, \DD)$, denote by $G(\DD_1, \DD)$, forms a group which acts on $B$ by left or right multiplication, and by conjugation. In particular, 
     $$
     \GD = \{\gm \in B \mid d_+(\gm) =1\}\,.
     $$
\end{propdef}
\begin{proof}
It follows directly from Proposition \ref{p:group} (i), (ii) and (iii) that $\HH(\DD_1, \DD):=(g_1 \KD g_1, B)$ is a standard table algebra with the degree map $d_+$. The remaining statements of the corollary are direct consequences of \cite[Prop. 3.18 and 3.19]{Bl}. We provide a proof here for the sake of completeness. \\
(i) $\Rightarrow$ (ii) and (iii): Let $\gm \in B$ be linear. Then $\gm \gmd = \ld g_1$ for some $\ld \in \BRp$. Then, by Proposition \ref{p:group} (iii) \eqref{eq:prod0}, $\ld = \frac{d_+(\MM)}{d_+(\DD_1)}$. On the other hand, we have
$$
\frac{d_+(\MM)^2}{d_+(\DD_1)^2} = d_+(\gm \gmd) = \ld d_+(g_1) = \ld = \frac{d_+(\MM)}{d_+(\DD_1)}\,.
$$
Since $\frac{d_+(\MM)}{d_+(\DD_1)} > 0$, $\ld = \frac{d_+(\MM)}{d_+(\DD_1)}=1$. Thus, $d_+(\gm)=1$ and $\gm\gmd=g_1$. \\
(ii) $\Rightarrow$ (i) and (iii): Since $\HH(\DD_1, \DD)$ is standard, if $d_+(\gm)=1$, then $C_{\MM, \MM^*}^{\DD_1} =1$. By Proposition \ref{p:group} (iii) \eqref{eq:prod0}, we have 
$$
1=d_+(\gm \gmd) =  d_+(g_1) +\sum_{\MM \in \I\setminus\{\DD_1\}} C_{\MM, \MM^*}^{\NN} d_+(\gn) \ge 1
$$
and hence $C_{\MM, \MM^*}^{\NN} =0$ for any $\NN \ne \DD_1$. Thus, $\gm \gmd = C_{\MM, \MM^*}^{\DD_1} g_1 = g_1$.

We continue to show that $G(\DD_1, \DD)$ is a group acts on $B$ by the left or right multiplication and conjugation. First, we  claim that $\gm \gn \in B$ for $\gm \in G(\DD_1, \DD)$ and $\gn \in B.$ By Proposition \ref{p:group} (iii) that
$$
\gm \gn = \sum_{\LL \in \I}\ld_\LL \gl \quad\text{where }\quad \ld_\LL \ge 0\,.
$$
By the preceding equivalences, we have
$$
\gn = \gmd\gm\gn = \sum_{\LL \in \I} C_{\MM, \NN}^\LL \gmd\gl\,.
$$
Every $\gmd\gl$ is a nonnegative linear combination of the elements in $B$. Therefore, $\gmd\gl$ must be a positive scalar multiple of $\gn$ whenever  $C_{\MM, \NN}^\LL > 0$. Since $\gmd$ is invertible,  the multiplication by $\gmd$ cannot send two distinct basis elements to proportional vectors. Therefore, $\gm \gn = \ld_\LL \gl$ for some $\LL \in \I$. On the other hand,
$$
\ld_\LL^2 [\gl, \gl] = [\gm \gn, \gm \gn] = [\gm^* \gm \gn,  \gn] = [\gn,  \gn]
$$
which implies $\ld_\LL^2 d_+(\LL) = d_+(\MM)$ or  $\ld_\LL^2 d_+(\gl) = d_+(\gm)$. Thus, $\ld_\LL^2 = \ld_\LL$. Since $\ld_\LL > 0$, $\ld_\LL=1$. Therefore, $\gm \gm = \gl$ and this proves the first claim. 

Second, we claim that $G(\DD_1, \DD)$ is closed under multiplication, which implies $G(\DD_1, \DD)$ is group. For any $\gm, \gn \in \GD$,
$\gm\gn = \gl$ for some $\LL \in B$ by the first claim, and $d_(\gl) = d_+(\gm \gn) =d_+(\gm) d_+(\gn) = 1$. Thus, $\gl \in \GD$, and so $\GD$ is closed under multiplication.

Now, the first claim implies $\GD$ acts on $B$ from the left by multiplication. Using the algebra antiautomorphism $*$, $B$ is closed under the $\GD$-action from the right. The closure of $B$ under the conjugation $\GD$-action follows directly from the left and the right $\GD$-actions on $B$.
\end{proof}

\begin{defn}\label{d:pseudocentral}
A fusion subcategory $\DD_1$ of a fusion category $\DD$ is called \emph{pseudocentral} if the idempotent 
$$
g_{\DD_1} = \frac{1}{d_+(\DD_1)} \sum_{X \in \irr(\DD_1)} d_+(X^*) X \in Z(\KD)\,.
$$
\end{defn}
\begin{rmk} \begin{enumerate}
\item  Let $H$ be a subgroup of a finite group $G$, $\w$ a 3-cocycle of $G$ and $\w_H$ its restriction on $H$. 
Then the fusion subcategory $\vec_H^{\w_H}$ is pseudocentral in $\vec_G^{\w}$ if and only if $H$ is a normal in $G$. 
\item Recall from \cite[Def. 2.4]{DMNO} that subcategory $\DD_1$ of $\DD$ is said to have a \emph{central structure} if the  inclusion functor $i : \DD_1 \to \DD$  factors through the forgetful functor $U:\cZ(\DD) \to \DD$, i.e. $i = U \tilde i$ for some tensor functor $\tilde i: \DD_1 \to \cZ(\DD)$. It obvious that $\DD_1$ is pseudocentral if it admits a central structure. \qedhere
\end{enumerate}
\end{rmk}
We can characterize $\GD$ if $\DD_1$ is a pseudocentral fusion subcategory of $\DD$ in terms of its objects.
\begin{thm}  \label{p:invertible}
Let $\DD_1$ be a pseudocentral fusion subcategory of any fusion category $\DD$, and let $\MM$ be an indecomposable left $\DD_1$-module subcategory of $\DD$. Then the following are equivalent:
\begin{enumerate}
    \item $\MM\otimes \MM^* \subseteq \DD_1$;
    \item $M \o M^* \in \DD_1$ for some $M \in \irr(\MM)$;
    \item $g_{\MM}$ is invertible.
\end{enumerate}
\end{thm}
\begin{proof} (i) $\Rightarrow$ (ii) is obviously true.\\
(ii) $\Rightarrow$ (iii): Since  $M\o M^* \in \DD_1$ for some $M \in \MM$, $C^\LL_{\MM,\MM^*}=0$ for all $\LL \ne \DD_1$ by Proposition \ref{p:group} (iv) \eqref{eq:prod1}, and so $\gm$ is linear. It follows from Proposition \ref{pd:table-alg} that  $\gm$ is invertible.\\
(iii) $\Rightarrow$ (i):  If there exist $M \in \MM$ and $N \in \MM^*$ such that $M \o N \not\in \DD_1$, then it follows from Proposition \ref{p:group} that
$$
g_1   \quad\text{and}\quad \gm \gmd=\frac{1}{d_+(M)d_+(N)} g_1 MN
$$ 
are linearly independent. This means $\gm$ is not linear, and so $\gm$ is not invertible by Proposition \ref{pd:table-alg}.   
\end{proof}

 We can extend the notion of group extensions of a fusion categories to table algebra extensions as follows. 
\begin{defn} \label{d:table-alg-extension}
Let $(\HH, B)$ be a table algebra. A fusion category $\DD$ is called an $(\HH, B)$-extension of a fusion category $\DD_1$ if there exists a table algebra homomorphism $\varphi: (K(\DD), \iD) \to (\HH, B)$. The extension is called \emph{faithful} if $\varphi$ is surjective.
\end{defn}

\begin{cor}\label{c:max-G-extension} Let $\DD$ be a fusion category, and $\DD_1$ a  fusion subcategory of $\DD$ such that $g_1 \in Z(\KD)$.  Then:
\begin{enumerate}
 \item The assignment $\varphi :X \mapsto g_1 X$ defines  an epimorphism of table algebras from $(\KD, \iD)$ to $\HH(\DD_1, \DD)$. In particular, $\DD$ is a faithful $\HH(\DD_1, \DD)$-extension of $\DD_1$. 
\item The full subcategory $E(\DD_1, \DD)$ linearly generated by
$$ 
\{X \in \iD \mid \varphi(X) = \ld_g g \text{ for some $g \in G$ and $\ld_g \in \BRp$} \}
$$
is a $G$-extension of $\DD_1$, where $G=\GD$, and $\ED$ is the largest faithful group extension of $\DD_1$ in $\DD$.
\end{enumerate}
\end{cor}
\begin{proof} (i) Since $g_1 \in Z(\KD)$, $g_1 \KD g_1 = g_1 \KD$.  It is obvious that map $\varphi(X)=g_1X$ for $X \in \iD$ defines an algebra epimorphism from $K(\DD)$ to $g_1K(\DD).$  For any $X \in \iD$, $X \in \irr(\MM)$ for some indecomposable $\DD_1$-bimodule $\MM$ in $\DD$. Thus, 
\begin{equation}
\varphi(X) = g_1 X = \frac{d_+(X) d_+(\DD_1)}{d_+(\MM)} \gm  = \frac{d_+(X)}{d_+(\gm)} \gm
\end{equation}
by Proposition \ref{p:group} (i) and (iii). Therefore, $\varphi$ is an epimorphism of table algebras.\\
(ii) Let $\I$ denote the set of indecomposable $\DD_1$-bimodule in $\DD$, and $\I_G =\{\MM \in \I\mid \gm \in G\}$.  Then, $\MM \mapsto \gm$ defines a bijection from $\I_G$ onto $G$. We claim that  
$$
\ED = \bigoplus_{\gm \in G} \MM\,.
$$
For any $\gm \in G$ and $X \in \irr(\MM)$,  $X \in \ED$ since
$$
\vp(X) =\frac{d_+(X)}{d_+(\gm)} \gm = d_+(X) \gm \,.
$$ 
Conversely, if $X$ is a simple object in $\ED$, then $g_1 X = \ld_\MM \gm$ for some $\gm \in G$ and $\ld_\MM > 0$. This implies $X \in \MM$, and we have proved the claim.

We continue to show that $\ED$ is a $G$-extension. Since $g_1 X = d_+(X) g_1$ for all $X \in \DD_1$, $\DD_1$ is a subcategory of $\ED$. Suppose $\gm, \gn \in G$, $X \in \MM$ and $Y \in \NN$. Then
$$
\vp(XY) =g_1XY = d_+(X) d_+(Y) \gm\gn = d_+(X) d_+(Y) \gl
$$
for some $\LL \in \I_G$ by Proposition \ref{pd:table-alg}.  Therefore, $\ED$ is a $G$-extension of $\DD_1$.

Suppose $\EE$ is a fusion subcategory of $\DD$ such that 
$$
\EE = \bigoplus_{h \in H} \EE_h
$$
is a faithful $H$-extension of $\DD_1 = \EE_1$ for some finite group $H$. Then $\EE_h$ is an indecomposable $\DD_1$-bimodule in $\DD$, and $d_+(\EE_h) = d_+(\DD_1)$ for all $h \in H$ by \cite[Thm. 3.5.2]{EGNO}. Let $g_h:=g_{\EE_h}$. Then $d_+(g_h) = d_+(\EE_h)/d_+(\DD_1) = 1$ and so $g_h \in \GD$. Obviously, $h \mapsto g_h$ defines a group monomorphism from $H$ to $G$. Therefore,
 $$
\EE = \bigoplus_{g_h \in G} \EE_h
$$
is a subgroup extension of $\DD_1$ in $\ED$. 
\end{proof}

\begin{rmk}  Theorem \ref{p:invertible} holds when $\DD=\CA$ and $\DD_1=\CAo$ where $A$ is condensable algebra in a modular fusion category $\CC$ as $\CAo$ admits a  central structure.
\end{rmk}

\section{Canonical basis of \texorpdfstring{$\CC(A,A)$}{}}
Let $A$ be a condensable algebra in a modular fusion category $\CC$. In this section, we continue the discussion of Section 6 for the subcategory $\CAo$  of $\CA$ to obtain a canonical basis $\ol{B}$ of $e_\1 \KA$  using the pivotal dimension character $d_A$ for $\KA$. This basis $\ol B$ is also indexed by the indecomposable $\CAo$-modules in $\CA$, the pair $(e_\1 \KA, \ol B)$ is a reality-based algebra, denoted $\oHA$. There is a bijection $\gm \mapsto \ogm$ from the distinguished basis $B$ of the table algebra $\HA$ to $\ol B$, and its restriction on $\GA$ is a group homomorphism.  The invertible elements in the table algebra $\HA$ form a group isomorphic to $\Auta(A)$ and the other elements of $\ol B$ give the noninvertible symmetries of $A$ via the canonical map $\rho$. 
We also identify the maximal group extension arising from categorical orbifold constructions.

Let $\DD$ be a pivotal fusion category,  $\DD_1$ a fusion subcategory of $\DD$  and let $\I$ denote the set of indecomposable $\DD_1$-bimodules in $\DD$. Recall from Section 6 that $\HD$ is a table algebra with the distinguished basis $B$ of   $\HD$ is given by 
$B=\{\gm \mid \MM \in \I\}$, where 
$$
\gm = \frac{1}{d_+(\MM)}\sum_{X \in \irr(\MM)} d_+(X) X\,.
$$
However, $\DD$ has a spherical structure which yields the categorical dimension character $d$ on $\KD$. One can define for each $\MM \in \I$ the element
$$
\ogm = \frac{1}{\dim(\MM)}\sum_{X \in \irr(\MM)} d(X) X
$$
where $\dim(\MM) = \sum_{X \in \irr(\MM)} d(X)^2$. Note that $\ol g_{\DD_1} = e_\1$, and we also denote it by $\og_1$. By using the same proof as in Proposition \ref{p:group} except (i), one could show $\og_1 \KD \og_1$ is an RBA which is analogous to Proposition \ref{pd:table-alg}.
\begin{propdef} \label{d:RBA} Let $\DD$ be a spherical fusion category,  $\DD_1$ a fusion subcategory of $\DD$  and let $\I$ denote the set of indecomposable $\DD_1$-bimodules in $\DD$. Then $e_\1 \KD e_\1$ is an RBA with the distinguished basis 
$$
\ol B=\{\ogm \mid \MM \in \I\}
$$
The degree map of this RBA is the restriction the categorical dimension character $d: \KD \to \BR$ and 
$$
d(\ogm)=d(\MM)/d(\DD_1) = C_{\MM, \MM^*}^{\DD_1}  
$$
for $\MM \in \I$, i.e., $d$ is standard. We will denote this standard RBA $(\og_1 \KD \og_1, \ol{B})$ by $\oHD$. Moreover, the map $\ol {\cdot} : \gm  \mapsto \ogm$ for $\MM \in \I$ defines a linear isomorphism from $\HD$ onto $\oHD$.
\end{propdef} 
\begin{proof}
By the remark preceding the statement, we need to provide another proof for $\og_1 \KM \og_1$ being 1-dimensional for any $\MM \in \I$ as $d(X)$ could be negative. This amount to the projection $P: \KM \to \KM$ defined by $P(M) = \og_1 M \og_1$ is a projection of rank 1. Since $\MM$, $\DD_1$ are subcategories of $\DD$, the spherical structure of $\DD$ induces a \emph{module trace} on $\MM$ (cf. \cite[Defn. 3.7]{Schau}). Now it follows from \cite[Props. 5.6 and 5.7]{Schau} that the matrix $Q:=([N, \og_1 M\og_1])_{M,N}$ indexed by $\irr(\MM)$ has rank 1. Since $Q$ is the matrix of $P$ relative to the basis $\irr(\MM)$, $P$ has rank 1. Following the same as the argument in Proposition \ref{p:group} (i), $Y\ogm = d(Y) \ogm = \ogm Y$ for $Y \in \DD_1$, and so $\og_1 \ogm \og_1= \ogm$. Thus, $\og_1 \KM \og_1$ is spanned by $\ogm$. 
\end{proof}

Here is our second main theorem in this paper: 
\begin{thm}\label{t:main2}
    Let $A$ be a condensable algebra in a modular fusion category $\CC$. 
    \begin{enumerate} 
    \item The algebra $\CC(A,A)$ is an RAR isomorphic to
    $\oHA$ with the distinguished basis $$\rho(\ol B) = \{\rho(\ogm)|\MM \in\I\}\,.$$ 
    For any $M \in \irr(\MM)$, 
    \begin{equation} \label{eq:rho_M}
    \frac{\dim(\MM)}{d_A(M)\dim(\CAo)} \rho(M) = \rho(\ogm)\,.
    \end{equation}
    In particular, for any $M,N \in \irr(\MM)$, we have
    $$
\frac{1}{d_A(M)} \rho(M) = \frac{1}{d_A(N)} \rho(N)\,.
    $$
    \item  $\Auta(A) \cong G(\CAo,\CA)$, and
    \begin{equation} \label{eq:AutA}
    \Auta(A)=\{\rho(\ogm)|d_+(\gm)=1 \text{ and }\MM \in \I \} \,.
    \end{equation}
    \item $\ECA$ is a twisted $A$-module subcategory given
    $$
    \ECA = \bigoplus_{g \in \Auta(A)} \CA(g)\,.
    $$    
    \end{enumerate}
\end{thm}
\begin{proof}
(i)   Since the fusion subcategory $\CAo$ admits a central structure, ${\ol g_1} = e_\1 \in Z(\KA)$ and $e_\1 \KA e_\1 = e_\1 \KA$.  Recall that the canonical algebra epimorphism $\rho: K(\CA) \to \CC(A,A)$ preserves duality and $e_\1 \KA \stackrel{\rho}{\cong} \CC(A,A)$ by \cite{DNRX2025} or Theorem \ref{t:Galois-fusion}. The image of the distinguished basis $\ol B$ of $e_\1 \KA$ under $\rho$ is a basis of $\CC(A,A)$. Now, the first assertion of (i) follows.

By Proposition \ref{d:RBA}, for any indecomposable $\CAo$-bimodule $\MM$ in $\CA$, $e_\1 K(\MM)$ is spanned by $\ogm$. Thus, for any $M \in \irr(\MM)$, $e_\1 M = \ld \ogm$ for some scalar $\ld \in \BR$. Applying $d_A$ to the equality, we find $d_A(M) = \ld \frac{\dim(\MM)}{\dim(\CAo)}$. By Theorem \ref{t:SW-fusion} (iv),  we find \eqref{eq:rho_M}.

(ii)  By Corollary \ref{c:max-G-extension}, $E(\CAo, \CA)$ is a $G$-extension of $\CAo$, where $G=G(\CAo, \CA)$. It follows from Theorem \ref{p5.17} that $g \mapsto \ol g$ defines an injective map from $G$ to $e_\1 \KA$ that preserve multiplications, and $\tau(g) =\rho(\og)$ defines  a group monomorphism from $G$ to $\Auta(A)$. We may simply identify $G$ as a subgroup of $\Auta(A)$ under $\tau$. By Theorem \ref{t2},
$$
\BB = \bigoplus_{h \in \Auta(A)} \CA(h)
$$
is   $(\CC_{A'}^0)_A$ is an $\Auta(A)$-extension of $\CAo$ in $\CA$, where $A'=A^{\Auta(A)}$, and hence $\Auta(A)$ is isomorphic to a subgroup of $G$ by Corollary \ref{c:max-G-extension}. Therefore, $G = \Auta(A)$. Since the elements of $\Auta(A)$ are of the form $\rho(\ogm)$ where $\gm \in G$, the equality \eqref{eq:AutA} follows from Proposition \ref{pd:table-alg}.\\
(iii) follows directly from (ii) and Corollary \ref{c:max-G-extension}.
\end{proof}
\begin{rmk} Under the conditions of Theorem \ref{t:main2}, if $d_A \ne d_+$, then $g_1 \og_1 = 0$. Since $\gm = g_1 \gm$ for any indecomposable $\CAo$-bimodule $\MM$ in $\CA$, $\rho(\gm) = \rho(\og_1 \gm) = 0$ by Theorem \ref{t:SW-fusion} (iv) as $e_\1 =\og_1$.
\end{rmk}

\begin{cor}
    Let $\DD$ be a spherical fusion category,  $I: \DD \to Z(\DD)$ the right adjoint of the forgetful functor $F:  Z(\DD) \to \DD$, and $A = I(\1_\DD)$. Then $\Aut_{Z(\DD)}^{\alg}(A)$ is isomorphic to the group of invertible objects of $\DD$. 
\end{cor}
\begin{proof}
    Since $A$ is a Lagrangian algebra in $\CC=Z(\DD)$, $\CAo = \Vec$ and $\CA$ is equivalent  to $\DD$ as spherical fusion categories. Therefore,  $E(\CAo, \CA)  \stackrel{\o}{\simeq} \DD_\pt$ and so $G(\Vec, \DD)$ is the group of invertible objects. The statement now follows immediately from Theorem \ref{t:main2} (ii).
\end{proof}

Recall the definition of braided $G$-crossed categories from \cite[8.24]{EGNO} or \cite[Def. 4.41]{DGNOI}. It has been stated in \cite{Ki3} that $G$-extensions of $\CAo$ by $g$-twisted $A$-modules for $g \in G$ are braided $G$-crossed extensions, and was proved in \cite{McRae2021}; here we provide a shorter proof for this result in our setting.  
\begin{cor} \label{c:bgc}
 Let $A$ be a condensable algebra in a modular fusion category $\CC$. If a fusion subcategory $\BB$ of $\CA$ is a $G$-extension for some finite group $G$, then $\BB$ is a braided $G$-crossed category. 
\end{cor}
\begin{proof}
   By Proposition \ref{p5.17}, $G$ is isomorphic to a subgroup of $\Auta(A)$. One can identify $G$ as a subgroup $\Auta(A)$ so that 
   $$
   \BB_g = \CA(g)\,.
   $$
   Let $B = A^G$. Then, by Proposition \ref{p5.17}, we have
   $$
   \BB = (\CBo)_A = \bigoplus_{g \in G} \CA(g)\,.
   $$
   Without losing generality, one may assume both $\CC_A$ and $\CC$ are strict pivotal. For $g \in G$, define the functor $T_g: \BB \to \BB$ by $T_g(X) = X$ in $\CC$ with the $A$-module structure $\mu_{T_g(X)} = \mu_X(\id_X \o g^{-1})$ and $T_g(f)=f$ for any morphism $f$ in $\BB$. Moreover, $g^{-1} : A \stackrel{\sim}{\to} T_g(A)$ in $\CA$, and  $T_g(X) \o_A T_g(Y) = T_g(X \o_A Y)$. To see this, one may recall that for any $X,Y \in \CA$, 
   $X \o_A Y$ is the image of the idempotent $P_{X,Y}: X \o Y \to X \o Y$ of $A$-modules, which has the diagrammatic presentation (cf. \cite[FIG. 14]{KO02})
\[
P_{X,Y} = \frac{1}{d(A)}\begin{tikzpicture} [baseline=-4ex]
\node (X) at (0,0) [above] { $\scriptstyle X $}; 
\node (Y) at (.8,0) [above] { $\scriptstyle Y $}; 
\begin{knot}[
end tolerance =1pt
clip width=5pt,
clip radius=5pt,
]
\strand[thick]  (X) to [out=270, in=90] +(0, -1.5);
\strand[thick]  (Y) to [out=270, in=90] +(0,-1.5);
\strand[densely dotted, thick] (0,-1)
to [out=45, in=90] (1.3,-.5)
to [out=-90, in=45] (.8,-1);
\end{knot}
\end{tikzpicture}\quad .
\]   
Since $g^* = g^{-1}$, we find
\[P_{T_g(X), T_g(Y)} =  \frac{1}{d(A)}\begin{tikzpicture} [baseline=-4ex]
\node at (.4,-.6) {$\scriptstyle g^{-1}$};
\node (g) at (.4,-.6) [draw, minimum width=15pt, minimum height=11pt,thick] {};
\node  at (1.2,-.6) {$\scriptstyle g^{-1}$};
\node (h) at (1.2,-.6) [draw, minimum width=15pt, minimum height=11pt,thick] {};
\node (X) at (0,0) [above] { $\scriptstyle X $}; 
\node (Y) at (.8,0) [above] { $\scriptstyle Y $}; 
\begin{knot}[
end tolerance =1pt
clip width=5pt,
clip radius=5pt,
]
\strand[thick]  (X) to [out=270, in=90] +(0, -1.5);
\strand[thick]  (Y) to [out=270, in=90] +(0,-1.5);
\strand[densely dotted, thick] (g.90)
to [out=80, in=100] (h.90);
\strand[densely dotted, thick] (h.270)
to [out=270, in=40]  (.8,-1.2);
\strand[densely dotted, thick] (g.270)
to [out=270, in=40]  (0,-1.2);
\end{knot}
\end{tikzpicture}
=\frac{1}{d(A)}\begin{tikzpicture} [baseline=-4ex]
\node (X) at (0,0) [above] { $\scriptstyle X $}; 
\node (Y) at (.8,0) [above] { $\scriptstyle Y $}; 
\begin{knot}[
end tolerance =1pt
clip width=5pt,
clip radius=5pt,
]
\strand[thick]  (X) to [out=270, in=90] +(0, -1.5);
\strand[thick]  (Y) to [out=270, in=90] +(0,-1.5);
\strand[densely dotted, thick] (0,-1)
to [out=45, in=90] (1.3,-.5)
to [out=-90, in=45] (.8,-1);
\end{knot}
\end{tikzpicture}
=P_{X,Y}
\]  
which implies $T_g \in \Aut_\o(\BB)$.   
   The functor $G \to \Fun(\BB, \BB)$, $g \mapsto T_g$, is obviously a strict monoidal functor by definition.  It is straightforward to verify that if $X \in \CA(h)$ then   
    $T_g(X)$ in $\CA(ghg^{-1})$.

    For any $X \in \CA(g)$ and $Y \in \CA$, the $G$-braiding $c_{X,Y}: X \o_A Y \to T_g(Y) \o_A X$ is inherited from $R_{Y,X}^{-1}$ and  is completely determined by the composition $R_{Y,X}^{-1} P_{X,Y}$, and $d(A) R_{Y,X}^{-1} P_{X,Y}$ 
\[
=\! \begin{tikzpicture} [baseline=-5ex]
\node (X) at (0,0) [above] { $\scriptstyle X $}; 
\node (Y) at (.6,0) [above] { $\scriptstyle Y $}; 
\begin{knot}[
end tolerance =1pt
clip width=10pt,
clip radius=10pt,
]
\strand[thick]  (X) to [out=270, in=90] +(0, -1.2)
to [out=270, in=90] +(0.6,-2);
\strand[thick]  (Y) to [out=270, in=90] +(0,-1.2)
to [out=270, in=90] +(-.6,-2);
\strand[densely dotted, thick] (0,-.8)
to [out=45, in=90] (1.1,-.3)
to [out=-90, in=45] (.6,-.8);
\flipcrossings{1}
\end{knot}
\end{tikzpicture}
\!=\!\begin{tikzpicture} [baseline=-5ex]
\node (X) at (0,0) [above] { $\scriptstyle X $}; 
\node (Y) at (.6,0) [above] { $\scriptstyle Y $}; 
\begin{knot}[
consider self intersections,
end tolerance =1pt
clip width=9pt,
clip radius=9pt,
]
\strand[thick]  (X) to [out=270, in=90] +(0, -1.2)
to [out=270, in=90] +(0.6,-2);
\strand[thick]  (Y) to [out=270, in=90] +(0,-1.2)
to [out=270, in=90] +(-.6,-2);
\strand[densely dotted, thick] (.6,-.6)
to [out=60, in=270] (1.4,-.1)
to [out=90, in=90] (1,-.3);
\strand[densely dotted, thick] (1,-.45)
to [out=270, in=60] (.5,-1.55);
\flipcrossings{1}
\end{knot}
\end{tikzpicture}
\!=\!\begin{tikzpicture} [baseline=-5ex]
\node (X) at (0,0) [above] { $\scriptstyle X $}; 
\node (Y) at (.6,0) [above] { $\scriptstyle Y $}; 
\begin{knot}[
consider self intersections,
end tolerance =1pt
clip width=9pt,
clip radius=9pt,
]
\strand[thick]  (X) to [out=270, in=90] +(0, -1.2)
to [out=270, in=90] +(0.6,-2);
\strand[thick]  (Y) to [out=270, in=90] +(0,-1.2)
to [out=270, in=90] +(-.6,-2);
\strand[densely dotted, thick] (.6,-.6)
to [out=60, in=180] (1,-.1)
to [out=0, in=60] (.5,-1.55);
\flipcrossings{1}
\end{knot}
\end{tikzpicture}
\hspace{-10pt}=\!\begin{tikzpicture} [baseline=-5ex]
\node (X) at (0,0) [above] { $\scriptstyle X $}; 
\node (Y) at (.6,0) [above] { $\scriptstyle Y $}; 
\begin{knot}[
consider self intersections,
end tolerance =1pt
clip width=10pt,
clip radius=10pt,
]
\strand[thick]  (X) to [out=270, in=90] +(.6, -.8)
to [out=270, in=90] +(0.6,-2);
\strand[thick]  (Y) to [out=270, in=90] +(-.6,-.8)
to [out=270, in=90] +(-.6,-2);
\strand[densely dotted, thick] (0,-1.4)
to [out=50, in=180] (.8,-.4)
to [out=0, in=60] (.6,-1.4);
\flipcrossings{1}
\flipcrossings{2}
\end{knot}
\end{tikzpicture}
\!\!=\!\!\begin{tikzpicture}[baseline=-5ex]
\node (X) at (0,0) [above] { $\scriptstyle X $}; 
\node (Y) at (.8,0) [above] { $\scriptstyle Y $}; 
\node (g) at (.4,-1) {$\scriptstyle g^{-1}$};
\node (g) at (.4,-1) [draw, minimum width=15pt,
minimum height=11pt,thick] {};
\begin{knot}[
consider self intersections,
end tolerance =1pt
clip width=10pt,
clip radius=10pt,
flip crossing/.list={1},
]
\strand[thick]  (X) to [out=270, in=90] +(.8, -.8)
to [out=270, in=90] +(.8,-2);
\strand[thick]  (Y) to [out=270, in=90] +(-.8,-.8)
to [out=270, in=90] +(-.8,-2);
\strand[densely dotted, thick] (0,-1.6)
to [out=50, in=270] (g);
\strand[densely dotted, thick] (g.90)
to [out=60, in=180] (.9,-.5)
to [out=0, in=45] (.8,-1.6);
\end{knot}
\end{tikzpicture} 
\!\!=\!\!\begin{tikzpicture}[baseline=-5ex]
\node (X) at (0,0) [above] { $\scriptstyle X $}; 
\node (Y) at (.8,0) [above] { $\scriptstyle Y $}; 
\node at (0,-1.75) { $\scriptstyle T_g(Y) $}; 
\node at (.8,-1.75) { $\scriptstyle X $}; 
\begin{knot}[
consider self intersections,
end tolerance =1pt
clip width=10pt,
clip radius=10pt,
]
\strand[thick]  (X) to [out=270, in=90] +(.8, -.8)
to [out=270, in=90] +(.8,-1.8);
\strand[thick]  (Y) to [out=270, in=90] +(-.8,-.8)
to [out=270, in=90] +(-.8,-1.8);
\strand[densely dotted, thick] (0,-1.3)
to [out=50, in=200] (1,-.6)
to [out=0, in=45] (.8,-1.3);
\flipcrossings{1}
\end{knot}
\end{tikzpicture} 
\,.
\] 
So, $c_{X,Y}: X \o_A Y \stackrel{\sim}{\to} T_g(X) \o_A Y$ with $c_{X,Y}^{-1}$ determined by $R_{Y,X}P_{T_g(Y),X}$.

To complete the proof of the statement, we need to verify the compatibility condition
\begin{equation} \label{eqL Gb1}
    c_{T_g(X), T_g(Y)} = c_{X,Y} = T_g(c_{X,Y})
\end{equation}
and two \emph{hexagon-like} equations  for any $X \in \CA(g)$ and $Y \in \CA(h)$:
\begin{eqnarray}
    \label{eq:hex1} c_{X, Y\o_A Z} & = & (T_g(Y) \o_A c_{X, Z}) (c_{X, Y} \o_A Z) \\
     \label{eq:hex2} c_{X\o_A Y, Z} & = & (c_{X, T_h(Z)} \o_A Y) (X \o_A c_{Y,Z})\,.
\end{eqnarray}

For any $g \in G$, $X\in \CA(h)$, $Y \in \CA$, $T_g(c_{X, Y}) = c_{X,Y}$ by definition, and $c_{T_g(X), T_g(Y)}$ is determined by $d(A) R_{T_g(Y),T_g(X)}^{-1} P_{T_g(X),T_g(Y)}$
$$
=\!\begin{tikzpicture} [baseline=-5ex]
\node (X) at (0,0){};
\node (0,0)  [above, xshift=-5pt] { $\scriptstyle T_g(X) $}; 
\node (Y) at (.6,0) {};
\node (.6,0) [above, xshift=20pt] { $\scriptstyle T_g(Y)$}; 
\begin{knot}[
end tolerance =1pt
clip width=10pt,
clip radius=10pt,
]
\strand[thick]  (X) to [out=270, in=90] +(0, -1.2)
to [out=270, in=90] +(0.6,-2);
\strand[thick]  (Y) to [out=270, in=90] +(0,-1.2)
to [out=270, in=90] +(-.6,-2);
\strand[densely dotted, thick] (0,-1)
to [out=45, in=90] (1.1,-.5)
to [out=-90, in=45] (.6,-1);
\flipcrossings{1}
\end{knot}\,.
\end{tikzpicture}
=
\begin{tikzpicture} [baseline=-5ex]
\node at (.4,-.6) {$\scriptstyle g^{-1}$};
\node (g) at (.4,-.6) [draw, minimum width=15pt, minimum height=11pt,thick] {};
\node  at (1.2,-.6) {$\scriptstyle g^{-1}$};
\node (h) at (1.2,-.6) [draw, minimum width=15pt, minimum height=11pt,thick] {};
\node (X) at (0,0) [above] { $\scriptstyle X $}; 
\node (Y) at (.8,0) [above] { $\scriptstyle Y $}; 
\begin{knot}[
end tolerance =1pt
clip width=9pt,
clip radius=9pt,
]
\strand[thick]  (X) to [out=270, in=90] +(0, -1.5)
to [out=270, in=90] +(.8,-2.1);
\strand[thick]  (Y) to [out=270, in=90] +(0,-1.5)
to [out=270, in=90] +(-.8,-2.1);
\strand[densely dotted, thick] (g.90)
to [out=80, in=100] (h.90);
\strand[densely dotted, thick] (h.270)
to [out=270, in=40]  (.8,-1.2);
\strand[densely dotted, thick] (g.270)
to [out=270, in=40]  (0,-1.2);
\flipcrossings{1}
\end{knot}
\end{tikzpicture}
=
\begin{tikzpicture} [baseline=-5ex]
\node (X) at (0,0) [above] { $\scriptstyle X $}; 
\node (Y) at (.6,0) [above] { $\scriptstyle Y $}; 
\begin{knot}[
end tolerance =1pt
clip width=10pt,
clip radius=10pt,
]
\strand[thick]  (X) to [out=270, in=90] +(0, -1.2)
to [out=270, in=90] +(0.6,-2);
\strand[thick]  (Y) to [out=270, in=90] +(0,-1.2)
to [out=270, in=90] +(-.6,-2);
\strand[densely dotted, thick] (0,-.8)
to [out=45, in=90] (1.1,-.3)
to [out=-90, in=45] (.6,-.8);
\flipcrossings{1}
\end{knot}
\end{tikzpicture} 
=d(A) R_{Y,X}^{-1} P_{X,Y} \,.
$$
Therefore, $c_{T_g(X), T_g(Y)} = c_{X,Y} = T_g(c_{X,Y})$.

\eqref{eq:hex1} follows from the following diagrammatical equalities
$$
\begin{tikzpicture} [baseline=-5ex]
\node (X) at (-.4,0) [above] { $\scriptstyle X $}; 
\node (Y) at (0,0) [above] { $\scriptstyle Y $}; 
\node (Z) at (.4,0) [above] { $\scriptstyle Z $}; 
\begin{knot}[
end tolerance =1pt
clip width=10pt,
clip radius=10pt,
]
\strand[thick]  (X) to [out=270, in=90] +(0, -1.2)
to [out=270, in=90] +(0.8,-2);
\strand[thick]  (Y) to [out=270, in=90] +(0, -1.2)
to [out=270, in=90] +(-0.4,-2);
\strand[thick]  (Z) to [out=270, in=90] +(0,-1.2)
to [out=270, in=90] +(-.4,-2);
\strand[densely dotted, thick] (0,-.5)
to [out=45, in=90] (.6,-.3)
to [out=-90, in=45] (.37,-.5);
\strand[densely dotted, thick] (-.4,-1)
to [out=45, in=90] (.6,-.8)
to [out=-90, in=45] (.37,-1);
\flipcrossings{1}
\flipcrossings{2}
\end{knot}
\end{tikzpicture} 
=
\begin{tikzpicture} [baseline=-5ex]
\node (X) at (-.4,0) [above] { $\scriptstyle X $}; 
\node (Y) at (0,0) [above] { $\scriptstyle Y $}; 
\node (Z) at (.4,0) [above] { $\scriptstyle Z $}; 
\begin{knot}[
end tolerance =1pt
clip width=10pt,
clip radius=10pt,
]
\strand[thick]  (X) to [out=270, in=90] +(0, -1.2)
to [out=270, in=90] +(0.8,-2);
\strand[thick]  (Y) to [out=270, in=90] +(0, -1.2)
to [out=270, in=90] +(-0.4,-2);
\strand[thick]  (Z) to [out=270, in=90] +(0,-1.2)
to [out=270, in=90] +(-.4,-2);
\strand[densely dotted, thick] (0,-.5)
to [out=45, in=90] (.6,-.3)
to [out=-90, in=45] (.37,-.5);
\strand[densely dotted, thick] (-.4,-1)
to [out=45, in=90] (.2,-.8)
to [out=-90, in=45] (0,-1);
\flipcrossings{1}
\flipcrossings{2}
\end{knot}
\end{tikzpicture} 
=
\begin{tikzpicture} [baseline=-5ex]
\node (X) at (-.4,0) [above] { $\scriptstyle X $}; 
\node (Y) at (0,0) [above] { $\scriptstyle Y $}; 
\node (Z) at (.4,0) [above] { $\scriptstyle Z $}; 
\begin{knot}[
end tolerance =1pt
clip width=10pt,
clip radius=10pt,
]
\strand[thick]  (X) to [out=270, in=90] +(0, -1.2)
to [out=270, in=90] +(0.8,-2);
\strand[thick]  (Y) to [out=270, in=90] +(0, -1.2)
to [out=270, in=90] +(-0.4,-2);
\strand[thick]  (Z) to [out=270, in=90] +(0,-1.2)
to [out=270, in=90] +(-.4,-2);
\strand[densely dotted, thick] (-0.4,-1)
to [out=45, in=90] (.6,-.7)
to [out=-90, in=45] (.37,-.9);
\strand[densely dotted, thick] (-.4,-.5)
to [out=45, in=90] (.2,-.3)
to [out=-90, in=45] (0,-.5);
\flipcrossings{1}
\flipcrossings{2}
\end{knot}
\end{tikzpicture} 
=
\begin{tikzpicture} [baseline=-5ex]
\node (X) at (-.4,0) [above] { $\scriptstyle X $}; 
\node (Y) at (0,0) [above] { $\scriptstyle Y $}; 
\node (Z) at (.4,0) [above] { $\scriptstyle Z $}; 
\begin{knot}[
end tolerance =1pt
clip width=10pt,
clip radius=10pt,
]
\strand[thick]  (X) to [out=270, in=90] +(0, -.6)
to [out=270, in=90] +(0.8,-2);
\strand[thick]  (Y) to [out=270, in=90] +(0, -.6)
to [out=270, in=90] +(-0.4,-2);
\strand[thick]  (Z) to [out=270, in=90] +(0,-1.2)
to [out=270, in=90] +(-.4,-2);
\strand[densely dotted, thick] (0,-1.1)
to [out=45, in=90] (.6,-.9)
to [out=-90, in=45] (.37,-1.1);
\strand[densely dotted, thick] (-.4,-.5)
to [out=45, in=90] (.2,-.3)
to [out=-90, in=45] (0,-.5);
\flipcrossings{1}
\flipcrossings{2}
\end{knot}
\end{tikzpicture} 
=
\begin{tikzpicture} [baseline=-5ex]
\node (X) at (-.4,0) [above] { $\scriptstyle X $}; 
\node at (-.4, -1.8) [xshift=-5pt] {$\scriptstyle T_g(Y) $};
\node (Y) at (0,0) [above] { $\scriptstyle Y $}; 
\node at (0, -1.8) {$\scriptstyle Z $};
\node (Z) at (.4,0) [above] { $\scriptstyle Z $}; 
\node at (.4, -1.8) [xshift=2pt] {$\scriptstyle X $};
\begin{knot}[
end tolerance =1pt
clip width=10pt,
clip radius=10pt,
]
\strand[thick]  (X) to [out=270, in=90] +(.4, -.7)
to [out=270, in=90] +(0.83,-1.8);
\strand[thick]  (Y) to [out=270, in=90] +(-.4, -.7)
to [out=270, in=90] +(-.4,-1.8);
\strand[thick]  (Z) to [out=270, in=90] +(0,-.7)
to [out=270, in=90] +(-.4,-1.8);
\strand[densely dotted, thick] (0,-.6)
to [out=45, in=90] (.6,-.4)
to [out=-90, in=45] (.37,-.6);
\strand[densely dotted, thick] (-.4,-1.5)
to [out=30, in=70] (.55,-1.4)
to [out=-90, in=70] (.39,-1.5);
\flipcrossings{1}
\end{knot}
\end{tikzpicture} \,.
$$
The equation \eqref{eq:hex2} follows from similar arguments.
\end{proof}
Now we can refine the statement of Corollary \ref{c:Pic1} by applying Corollary \ref{c:bgc} and \cite[Thm. 7.12]{ENO2010}
\begin{cor} \label{c:Pic2}
    Let $A$ be a condensable algebra in modular fusion category $\CC$. Then there is a morphism of 2-groups $c : \Auta(A) \to \underline{\underline{\rm Pic}}(\CAo) $. \qedhere
\end{cor}
\begin{cor}\label{c:regular}
    Let $G$ be a subgroup of  $\Auta(A)$, and $B = A^G$. Then $\Rep(G)$ is equivalent to an isotropic fusion subcategory $\DD$ of $\CBo$ and $A$ is a regular algebra in $\DD$, where $\irr(\DD) = \{X \in \irr(\CBo)\mid [A, X]_{\CBo} \ne 0\}$. 
\end{cor}
\begin{proof}
   By Theorem \ref{t2}, $G= \Aut^{\alg}_{\CBo}(A)$ and $A^G = B$ in $\CBo$. According to Corollary \ref{c:bgc}, $(\CBo)_A$ is a braided $G$-crossed category. Now, the result follows directly from \cite[Thms. 2.11 and 2.15]{Ki2002}.
\end{proof}
 The corollary also follows from \cite[Thm.4.44]{DGNOI} as $(\CBo)_A$ is braided $G$-crossed category. By using Corollary \ref{c:degree}, the result can be proved by other approaches  \cite[Theorem 9]{DM2021}, \cite[Theorem 6.1]{DNR2025}. 
 
 The following lemma should be well-known but we include here for the sake of completeness.
 \begin{lem} \label{l:aut}
     Let $G$ be a finite group, $\EE=\Rep(G)$, $H$ a subgroup of $G$, and let $N_G(H)$ denote the normalizer of $H$ in $G$. Then $N_G(H)/H \cong \Aut_\EE^{\alg}(\BC[G/H]^*)$ under the assignment $\ol n \mapsto \sigma_{\ol n}$ with $\sigma_{\ol n}(e_{\ol g}) = e_{\ol {gn^{-1}}}$ for $\ol g \in G/H$ where $e_{\ol g} \in \BC[G/H]^*$ is defined by $e_{\ol g}(\ol k) = \delta_{\ol g, \ol k}$. In particular, $G \cong \Aut_\EE^{\alg}(\BC[G]^*)$ and $(\BC[G]^*)^H \cong \BC[G/H]^*$ as algebras in  $\EE$.
 \end{lem}
 \begin{proof}
     Note that $\BC[G/H]^*$ is a commutative semisimple algebra with the set of primitive idempotents $\{e_{\ol g}\mid \ol g \in G/H\}$. Thus, for any permutation $ \sigma \in \Sym(G/H)$, the linear extension of the assignment $\hat \sigma : e_{\ol g} \mapsto e_{\sigma \ol g}$ defines an algebra automorphism of $\BC[G/H]^*$ and every algebra automorphism of $\BC[G/H]^*$ is of this form. So we can identify $\Sym(G/H)$ with $\Aut_{\Vec}^{\alg}(\BC[G/H]^*)$. Since the $G$-action on $\BC[G/H]^*$ is given by $ae_{\ol g} = e_{\ol{ag}}$ for any $a, g \in G$, $\Aut_\EE^{\alg}(\BC[G/H]^*)$ is the group of $G$-equivariant permutations $\Sym_G(G/H)$.
     For $\s \in \Sym_G(G/H)$, $\s$ is completely determined by $\s(\ol 1) = \ol n$ as $\s(\ol g) = g \s(\ol 1) = \ol {gn}$. Since $\ol n= \s(\ol h) =  \ol{hn}$ for all $h \in H$, $n \in N_G(H)$. Conversely, for any $\ol n \in N_G(H)/H$, $\s_{\ol n} (\ol g) = \ol {g n^{-1}}$ defines a $G$-equivariant permutation on $G/H$. Moreover, the assignment $N_G(H)/H \to \Sym_G(G/H)$, $\ol n \mapsto \s_{\ol n}$, defines a group isomorphism.

     Applying to the preceding isomorphism to the trivial subgroup $\{1\}$, we find  $G \cong \Aut_\EE^{\alg}(\BC[G]^*)$. For any subgroup $H$ of $G$,  the idempotent $e_{gH}:=\sum_{a \in gH} e_a$ is fixed by $\s_h$ for all $h \in H$. The set of idempotents $\{e_{gH} \mid gH \in G/H\}$ spanned an algebra isomorphic to $\BC[G/H]^*$ in $\EE$. Thus, $(\BC[G]^*)^H \cong \BC[G/H]^*$ as algebras in $\EE$ by comparing their dimensions.
 \end{proof}
 \begin{cor}\label{c:xu}
      Let $G$ be a subgroup of  $\Auta(A)$, $H$ a subgroup of $G$, and $B = A^G$. Then, $A^H$ is a condensable subalgebra of $A$ and $\Aut_{\CBo}^{\alg}(A^H) \cong N_G(H)/H$.
 \end{cor}
 \begin{proof}
     The first assertion follows directly from Proposition \ref{p:5.16} (ii). For the second assertion, it suffices to show that if $B=\1$, then $\Auta(A^H) \cong N_G(H)/H$. In this case, $A$ is the regular algebra of the Tannakian subcategory $\DD$ of $\CC$ by Corollary \ref{c:regular}, and $\DD$ is equivalent to $\Rep(G)$ as premodular categories. Let $\EE = \Rep(G)$. Then, 
     $$
     \Auta(A^H) \cong \Aut_\EE^{\alg}((\BC[G]^*)^H) \cong \Aut_\EE^{\alg}(\BC[G/H]^*)  
     $$
     which is isomorphic $N_G(H)/H$ by Lemma \ref{l:aut}.
 \end{proof}

\bibliographystyle{plain}

\end{document}

%% file: Macro.tex
\usepackage{amsmath, amsfonts, amssymb, amsthm, color, dsfont}
\usepackage{enumitem}
\usepackage[utf8]{inputenc}
\usepackage[colorlinks,linkcolor=blue,citecolor=blue]{hyperref}
\usepackage[all, knot, color]{xy}
\xyoption{arc}
\xyoption{line}
\usepackage{booktabs}
\UseCrayolaColors

\usepackage{tikz}
\usetikzlibrary{cd, calc, knots, shapes.geometric, intersections}
\tikzset{
  dotted stroke/.style={thick, densely dotted}
}

\usepackage{euscript}

\newcommand\eB           {\EuScript{B}}
\newcommand\eC           {\EuScript{C}}
\newcommand\eD           {\EuScript{D}}

\newcommand\eL          {\EuScript{L}}
\newcommand\eM          {\EuScript{M}}
\newcommand\eN         {\EuScript{N}}

\newcommand\cI           {\mathcal{I}}

\newcommand\cZ         {\mathcal{Z}}

\renewcommand{\i}{{\iota}}
\newcommand{\BB}{{\eB}}

\newcommand{\HH}{\mathcal{H}}

\newcommand{\CC}{{\eC}}
\newcommand{\LL}{{\eL}}
\newcommand{\NN}{{\eN}}

\newcommand{\DD}{{\eD}}

\newcommand{\EE}{\mathcal{E}}

\newcommand{\MM}{{\eM}}

\newcommand{\BC}{\mathbb{C}}
\newcommand{\BN}{\mathbb{N}}
\newcommand{\BQ}{\mathbb{Q}}
\newcommand{\BZ}{\mathbb{Z}}

\newcommand{\BR}{\mathbb{R}}

\newcommand{\s}{\sigma}
\renewcommand{\b}{\beta}

\newcommand{\KD}{{K(\DD)}}
\newcommand{\KA}{{K(\CA)}}

\newcommand{\ld}{\lambda}

\newcommand{\e}{\varepsilon}
\newcommand{\vp}{\varphi}

\newcommand{\ev}{\operatorname{ev}}
\newcommand{\coev}{\operatorname{coev}}

\newcommand{\id}{\operatorname{id}}
\newcommand{\Hom}{\operatorname{Hom}}

\newcommand{\End}{\operatorname{End}}
\newcommand{\gal}{\operatorname{Gal}}
\newcommand{\irr}{\operatorname{Irr}}

\newcommand{\rev}{\operatorname{rev}}
\newcommand{\tr}{\operatorname{tr}}
\newcommand{\Tr}{\operatorname{Tr}}

\newcommand{\Rep}{\operatorname{Rep}}

\renewcommand{\1}{{\mathds{1}}}

\renewcommand{\Vec}{\operatorname{Vec}}

\newcommand{\pt}{{\operatorname{pt}}}

\usepackage{thmtools}

\newtheorem{thm}{Theorem}[section]
\newtheorem{lem}[thm]{Lemma}
\newtheorem{prop}[thm]{Proposition}
\newtheorem{cor}[thm]{Corollary}

\theoremstyle{definition}
\newtheorem{defn}[thm]{Definition}

\declaretheorem[
  name=Remark,
  style=definition,
  sibling=thm,
  qed=\ensuremath{\blacksquare}
]{rmk}

\def\I{{\cI}}

\def\GQ{\gal(\ol\BQ)}

\def\O2{{\Omega^n_2}}

\def\Sym{\operatorname{Sym}}

\def\alg{\operatorname{alg}}

\def\a{{\alpha}}

\def\o{{\otimes}}
\def\oa{{\otimes\!_A}}
\def\ol{\overline}

\def\w{\omega}
\def\CAo{{\CC_A^{\hspace{1pt}0}}}
\def\CA{{\CC_A}}
\def\CB{{\CC_B}}
\def\CBo{{\CC_B^0}}

\def\iC{\irr(\CC)}
\def\iD{{\irr(\DD)}}
\def\iCA{{\irr(\CC_A)}}
\def\iCAo{{\irr(\CAo)}}
\def\iCBo{{\irr(\CBo)}}

\def\Aut{{\rm Aut}}
\def\Fun{{\operatorname{Fun}}}

\def\Auta{{\operatorname{Aut}_\CC^{\operatorname{alg}}}}
\def\AutC{{\operatorname{Aut}_{\CC}}}
\def\gm{{g_\MM}}
\def\gmd{{g_{\MM^*}}}
\def\gn{{g_\NN}}
\def\gl{{g_\LL}}
\def\BRp{{\BR_{> 0}}}
\def\BRn{{\BR_{\ge 0}}}
\def\GD{{G(\DD_1, \DD)}}
\def\ED{{E(\DD_1, \DD)}}
\def\HD{{\HH(\DD_1, \DD)}}
\def\oHD{{\ol{\HH}(\DD_1, \DD)}}
\def\ECA{{E(\CAo, \CA)}}
\def\oHA{{\ol{\HH}(\CAo, \CA)}}
\def\HA{{\HH(\CAo, \CA)}}
\def\GA{{G(\CAo, \CA)}}
\def\og{{\ol{g}}}
\def\ogm{{\ol{g}_\MM}}
\def\KM{{K(\MM)}}

\def\vec{{\operatorname{Vec}}}

\usepackage[normalem]{ulem}
\usepackage{relsize, stackengine, graphicx}
\stackMath

\setlist[enumerate, 1]{label=\rm{(\roman*)}, leftmargin=18pt, labelsep=3pt}
\setlist[enumerate, 2]{label=\rm{(\alph*)}, leftmargin=18pt, labelsep=3pt}
\newenvironment{enumeri}
{\begin{enumerate}[label=\rm{(\roman*)}, leftmargin=18pt, labelsep=3pt]}
{\end{enumerate}}